\documentclass[12pt,reqno]{amsart}

\usepackage[english]{babel}
\usepackage[
  letterpaper,
  top=3cm,
  bottom=3cm,
  left=3cm,
  right=3cm,
  marginparwidth=1.75cm
]{geometry}

\usepackage{amsmath,amssymb,amsfonts,amsthm,amscd,amsbsy,mathrsfs}

\usepackage{graphicx}
\graphicspath{{figs2/}}

\usepackage{booktabs}
\usepackage{multirow}
\usepackage{subfigure}
\usepackage{xcolor}
\usepackage{textcomp}
\usepackage[title]{appendix}
\usepackage{manyfoot}

\usepackage{algorithm}
\usepackage{algorithmicx}
\usepackage{algpseudocode}
\usepackage{listings}

\usepackage{orcidlink}

\newtheorem{theorem}{Theorem}

\theoremstyle{definition}

\theoremstyle{remark}

\usepackage[square,numbers]{natbib}

\begin{document}

\title[Equilateral Restricted Four-Body Problem]
{Nonlinear Stability, Resonances, and Singular Reduction in the Unequal-Mass Equilateral Restricted Four-Body Problem}

\author[M. \'Alvarez-Ram\'irez]{Martha \'Alvarez-Ram\'irez}
\address{
Departamento de Matem\'aticas\\
Universidad Aut\'onoma Metropolitana--Iztapalapa\\
Mexico City, Mexico\\
\orcidlink{0000-0001-9187-1757}
}
\email{mar@xanum.uam.mx}

\author[J. A. Zepeda-Ram\'irez]{J. Alejandro Zepeda-Ram\'irez}
\address{
Departamento de F\'isica\\
Universidad Aut\'onoma Metropolitana--Iztapalapa\\
Mexico City, Mexico\\
\orcidlink{0000-0002-8189-4792}
}
\email{zeral@xanum.uam.mx}

\begin{abstract}
We study the nonlinear stability of the elliptic equilibrium points $L_3$, $L_5$, and $L_6$
in the planar equilateral restricted four-body problem with unequal primary masses. Existing
studies of this problem have addressed either equal-mass or partially symmetric configurations,
or have treated local stability criteria and global dynamics separately; a systematic treatment
of the genuinely unequal-mass case across the full two-parameter mass plane has been missing.
The present work fills this gap by combining three tools that had not previously been used
together in this setting: high-order Birkhoff normal forms, systematic numerical continuation
over the admissible mass plane, and singular reduction of the resonant normal forms. This
combination yields, for the first time, a complete classification of nonlinear stability
throughout the admissible two-parameter mass plane, rather than at isolated mass values.
In nonresonant regions, stability is established through Arnold's theorem whenever the
quartic nondegeneracy condition is satisfied; along the associated degeneracy curves, the
normalization is extended to sixth order to resolve the inconclusive cases. The low-order
resonances of type $2$:$-1$  and $3$:$-1$ are analyzed with the classical
theorems of Alfriend, Markeev, and Meyer, providing the first explicit classification of the
$3$:$-1$ resonant dynamics in the unequal-mass problem. To connect these local
results to the global phase-space geometry, we perform singular reduction of the resonant
normal forms and obtain the reduced orbit spaces, revealing the saddle-center bifurcation
mechanisms responsible for the appearance and disappearance of families of periodic orbits.
Taken together, these results extend previous nonlinear stability analyses beyond symmetric
mass configurations and provide a geometric interpretation of the resonant dynamics in the
unequal-mass equilateral restricted four-body problem.
\end{abstract}



\maketitle

\section{Introduction}\label{sec0}
In this work, we consider the planar equilateral restricted four-body problem (in short, ERFBP). Three massive primaries move on uniform circular orbits in a common plane while maintaining a rigid equilateral (Lagrange-type) central configuration. A fourth, infinitesimal test particle moves in the same plane under the Newtonian gravitational field generated by the primaries; its influence on their motion is neglected.

Over the past decades, the ERFBP has attracted growing interest, with
numerous studies devoted to its equilibria, their stability, and the
associated global dynamics. Much of this literature has focused on the
equal-mass case or on the subclass with two equal masses, as these
symmetric settings greatly facilitate the analysis (see, e.g.,
\cite{Balta1,AlvarezBarrabes2014,zotos} and references therein). By
contrast, the genuinely unequal-mass problem has received comparatively
little attention, even though it exhibits a
much richer dynamical structure.
The lack of symmetry reductions, the presence of two independent mass
parameters, and the resulting higher-dimensional parameter space
significantly increase both the analytical and computational
complexity. In this regime, the system naturally exhibits
a richer resonance and degeneracy structure, which
plays a fundamental role in determining its nonlinear stability
properties.

Pioneering contributions to the unequal-mass ERFBP were made by
Pedersen~\cite{pede}, Simó~\cite{simo}, and later by
Leandro~\cite{leandro}, followed by the numerical investigations of
Baltagiannis and Papadakis~\cite{Balta1}, who mapped the distribution
of equilibrium points and showed that, although most are linearly
unstable, the elliptic equilibria $L_3$, $L_5$, and $L_6$ may become
linearly stable within certain ranges of the mass parameters.
Subsequent work examined nonlinear stability under specific symmetry
assumptions. In the two-equal-masses case,
Álvarez-Ramírez et al.~\cite{alvastuchi} established nonlinear
stability through the construction of the Birkhoff normal form, while
Zepeda-Ramírez et al.~\cite{zepAl} extended this analysis to selected
unequal-mass configurations (restricted to
$m_1=0.01$ and $m_1=0.99$) and identified representative stability
domains.
Parallel developments by Bardin and
Volkov~\cite{BardinVolkov2024} investigated the nonlinear stability of
the ERFBP using Hamiltonian normal forms and KAM theory,
providing detailed resonance curves and stability diagrams in the admissible
mass plane. Their analysis established rigorous stability and
instability criteria based on local normal forms, both in resonant and
nonresonant regimes.
These works demonstrate that linear stability alone is
insufficient to characterize the dynamics near the elliptic
equilibria, and that resonances together with higher-order terms play
a central role in the nonlinear stability of the ERFBP.

The present work departs from previous studies in its emphasis. Rather than focusing exclusively on the derivation of local stability criteria, we integrate high-order Birkhoff normal forms, systematic numerical continuation over the admissible two-parameter mass plane, and singular reduction into a single framework for the nonlinear stability analysis of the unequal-mass ERFBP. This integration yields a systematic classification of the nonlinear stability of the elliptic equilibria across the mass plane, together with a geometric interpretation of the associated resonant dynamics. The analysis combines three complementary components.

First, we compute the quartic Arnold determinant $D_4$ and the
sixth-order determinant $D_6$, thereby determining the nonlinear
stability boundaries of the elliptic equilibria in the
$(m_1,m_2)$ parameter plane. Second, we construct the low-order
resonance curves and apply the Alfriend-Markeev-Meyer  theorem~\cite{cabral} and the
Alfriend-Markeev theorem~\cite{MeyerOffin} to classify the
corresponding resonant dynamics. Finally, we perform singular reduction
of the resonant Hamiltonians to obtain the associated  orbit
spaces and describe the bifurcation mechanisms of the resulting
families of periodic orbits.

The first component combines perturbation theory with high-order normal
form computations. We construct Birkhoff normal forms around the
elliptic equilibria and apply Arnold's theorem to establish nonlinear
stability in nonresonant regions, verifying the nondegeneracy
conditions through systematic numerical exploration of the
$(m_1,m_2)$ parameter plane. When the quartic Arnold determinant
$D_4$ vanishes, the normalization is extended to sixth order, leading
to the determinant $D_6$, which provides the corresponding stability
criterion. Resonant cases are treated using the
Alfriend-Markeev-Meyer theorem~\cite{cabral} for the
$2$:$-1$ resonance and the Alfriend-Markeev
theorem~\cite{MeyerOffin} for the $3$:$-1$ resonance, providing
rigorous criteria for the classification of the dynamics near these
commensurabilities.

To complement the local analysis, the second part of the paper adopts a
global perspective based on singular reduction. In action-angle
variables, resonances manifest through specific linear combinations of
the angular coordinates, which generate an $\mathbb{S}^1$ rotational
symmetry. Reduction with respect to this symmetry yields orbit spaces
endowed with stratified symplectic structures, providing a natural
geometric framework for the analysis of resonant dynamics. Formally,
the orbit space $\mathbb{O}$ is obtained from the energy surface
\[N_h=\{z\in M : H(z)=h\},\]
by projecting each trajectory of the Hamiltonian flow to a single
point, thereby encoding the reduced dynamics, where $M$ denotes the
full symplectic phase space of the system. Local symplectic
coordinates, expressed in terms of polynomial invariants, are
introduced around both regular and singular points of $\mathbb{O}$,
offering a desingularized setting for systematic stability and
bifurcation analysis. This approach, originally developed in the works
of Arms, Cushman, and Gotay~\cite{arms1991universal}, and subsequently
applied to celestial mechanics by Meyer, Palacián, Yanguas, and
collaborators~\cite{MeyerOffin,palacian2015,palacian2018, FerrerPalacian2026},
provides a natural framework to describe how resonance curves and
degeneracy loci organize the dynamics. In our context, we perform singular
reduction of the resonant Hamiltonians to obtain the associated orbit spaces and describe the bifurcation mechanisms of the resulting families of periodic orbits.

From a scientific standpoint, the main contribution of this work is threefold:
it extends nonlinear stability results for the ERFBP beyond the symmetric
mass configurations previously treated in the literature to the full
two-parameter mass plane; it provides the first explicit classification of
the $3$:$-1$ resonant dynamics in the unequal-mass case; and it establishes,
via singular reduction, a direct geometric link between local normal-form
stability and the global bifurcation structure of the associated periodic
orbit families. The present analysis is restricted to the low-order
$2$:$-1$ and $3$:$-1$ resonances and to the planar problem; higher-order
resonances and the spatial ERFBP lie beyond its scope and are discussed
briefly in the concluding remarks.

Beyond this restriction in scope, two further methodological limitations should be
stated explicitly. First, the stability theorems of Arnold, Alfriend, Markeev, and Meyer employed here are local:
they characterize the dynamics in a neighborhood of the equilibrium
and do not, by themselves, describe the size of the region of stability or its dependence on
the mass parameters. Second, the classification relies on high-precision numerical evaluation
of the normal-form coefficients rather than closed-form expressions, since the equilibrium
coordinates of the unequal-mass ERFBP are not available analytically; although this evaluation
was systematically checked for consistency under changes in working precision, it does not
constitute a rigorous computer-assisted proof. These limitations delineate natural directions
for future work, discussed further in the concluding remarks.

The paper is organized as follows. Section~\ref{sec:back} recalls the
formulation of the ERFBP, the equilibrium points, and their linear
stability properties. Section~\ref{sec:birkhoffnorma} develops the
Birkhoff normal form and studies the nonlinear stability in
nonresonant and degenerate cases. Section~\ref{sec_nonstab} analyzes
the low-order $2\!:\!-1$ and $3\!:\!-1$ resonances. Finally,
Section~\ref{sec:singular} presents the singular reduction of the
resonant Hamiltonians together with the corresponding reduced orbit
spaces and bifurcation diagrams.

\section{Mathematical setting}\label{sec:back}
In this section, we briefly recall the formulation of the planar equilateral restricted
four-body problem (ERFBP) and summarize the properties of its equilibrium points
that are needed for the subsequent nonlinear stability analysis.

\subsection{Mathematical formulation of the ERFBP}
We consider the motion of a massless particle under the gravitational
influence of three massive bodies (the primaries), which move in a
Lagrange equilateral configuration, preserving their mutual distances
while revolving uniformly about their common center of mass.

The Lagrange equilateral configuration is a classical central
configuration of the three-body problem and therefore remains rigid
while rotating with constant angular velocity. Exploiting the
rotational, translational, and scaling invariance of the Newtonian
equations, we normalize the angular velocity to $\omega=1$, the mutual
distance between the primaries to one, and place the center of mass at
the origin of a uniformly rotating reference frame. Without loss of
generality, the configuration is oriented so that the primary $m_1$
lies on the positive $x$-axis. Under these assumptions, the primaries
remain fixed in the rotating frame, while only the infinitesimal
particle evolves dynamically.

The equations of motion of the infinitesimal particle in the rotating
frame are
\begin{equation}\label{four1}
\begin{split}
\ddot{x}-2\dot{y} &= \Omega_x,\\
\ddot{y}+2\dot{x} &= \Omega_y,
\end{split}
\end{equation}
where dots denote derivatives with respect to time $t$. The terms
$2\dot{y}$ and $2\dot{x}$ are the Coriolis forces associated with the rotating frame, whereas the derivatives of the effective
potential $\Omega$ include both the gravitational attraction of the primaries and the centrifugal
potential. Working in this rotating frame is what allows the primaries themselves to
appear as fixed points, so that equation \eqref{four1} fully describes the dynamics of the
fourth, infinitesimal body relative to a static background configuration.

The positions of the primaries are determined uniquely by the
constraints that they form an equilateral triangle of unit side length,
their center of mass is located at the origin, and the primary $m_1$
lies on the positive $x$-axis. These conditions are expressed as
\[
\begin{aligned}
&(x_1-x_2)^2+(y_1-y_2)^2=1,\qquad
(x_1-x_3)^2+(y_1-y_3)^2=1,\\
&(x_2-x_3)^2+(y_2-y_3)^2=1,\qquad
m_1x_1+m_2x_2+m_3x_3=0,\\
&m_1y_1+m_2y_2+m_3y_3=0,\qquad
y_1=0,
\end{aligned}
\]
which admit a unique solution for every admissible choice of the masses (see \cite{Moulton1900}).
Physically, this system fixes the shape of the equilateral triangle once the
mass parameters are given, up to the two-fold reflection ambiguity
corresponding to the two possible orientations of the triangle; the
condition $y_1=0$, together with $m_1$ lying on the positive $x$-axis,
selects one of these orientations, so that the coordinates of the
primaries are fixed unambiguously once this convention is adopted.
Consequently, the coordinates of the primaries are
\begin{equation}\label{eq:prim}
\begin{array}{lll}
&x_{1} = \dfrac{|K|\sqrt{m^{2}_{2}+m_{2}m_{3}+m^{2}_{3}}}{K},
\qquad  & y_{1}=0,\vspace{0.3cm}\\
&x_{2} =
-\dfrac{|K|\left[(m_{2}-m_{3})m_{3}+m_{1}(2m_{2}+m_{3})\right]}
{2K\sqrt{m^{2}_{2}+m_{2}m_{3}+m^{2}_{3}}},
\qquad
& y_{2}= \dfrac{\sqrt{3}}{2}
\dfrac{m_3}
{\sqrt{m^{2}_{2}+m_{2}m_{3}+m^{2}_{3}}},
\vspace{0.3cm}\\
& x_{3} =
-\dfrac{|K|}{2\sqrt{m^{2}_{2}+m_{2}m_{3}+m^{2}_{3}}},
\qquad
&y_{3}=
- \dfrac{\sqrt{3}}{2}
\dfrac{m_2}
{\sqrt{m^{2}_{2}+m_{2}m_{3}+m^{2}_{3}}},
\end{array}
\end{equation}
where
\[
K=m_{2}(m_{3}-m_{2})+m_{1}(m_{2}+2m_{3}),
\]
and the masses satisfy the normalization
\begin{equation}\label{eq_mass}
m_1+m_2+m_3=1.
\end{equation}
This normalization fixes the overall mass scale so that only the relative proportions among 
the three primaries matter dynamically, reducing the problem to the two independent parameters 
$(m_1,m_2)$ used throughout the paper.

The effective potential is
\[
\Omega(x,y)=
\frac{1}{2}(x^2+y^2)
+\frac{m_1}{r_1}
+\frac{m_2}{r_2}
+\frac{m_3}{r_3},
\]
where
\[
r_1=\sqrt{(x-x_1)^2+y^2},\qquad
r_2=\sqrt{(x-x_2)^2+(y-y_2)^2},\qquad
r_3=\sqrt{(x-x_3)^2+(y-y_3)^2}.
\]

Introducing the conjugate momenta
\[
p_x=\dot{x}-y,
\qquad
p_y=\dot{y}+x,
\]
the Hamiltonian takes the form
\begin{equation}\label{ham1}
H=
\frac{1}{2}(p_x^2+p_y^2)
+y p_x
-x p_y
-U(x,y),
\end{equation}
where
\[
U(x,y)=
\frac{m_1}{r_1}
+\frac{m_2}{r_2}
+\frac{m_3}{r_3}.
\]
Using \eqref{eq_mass}, we write
\[
m_3=1-m_1-m_2,
\]
so the problem depends on two independent mass parameters, $m_1$ and $m_2$.

Equations~\eqref{four1}, \eqref{eq:prim} and \eqref{eq_mass} reduce the problem of locating the primaries and
describing the motion of the infinitesimal particle to a system with only two free
parameters, once symmetry and scaling have been used. These two parameters, $(m_1,m_2)$, are the coordinates used
throughout the remainder of the article: every equilibrium point,
resonance curve, and stability boundary discussed below appears as a
curve or region in this same mass plane.

\subsection{Equilibria and linear stability analysis}
The equilibrium points are determined by the equations
\begin{equation}\label{lss1}
\Omega_x(x,y,m_1,m_2)=
\Omega_y(x,y,m_1,m_2)=0.
\end{equation}
Since these equations cannot be solved analytically, all computations were carried
out numerically using \textit{Mathematica} with variable-precision arithmetic,
typically between 40 and 60 digits. The computations were systematically verified
by increasing the working precision, and the qualitative conclusions were found to
be robust.

\begin{figure}[H]
\centering
\includegraphics[width=0.6\textwidth]{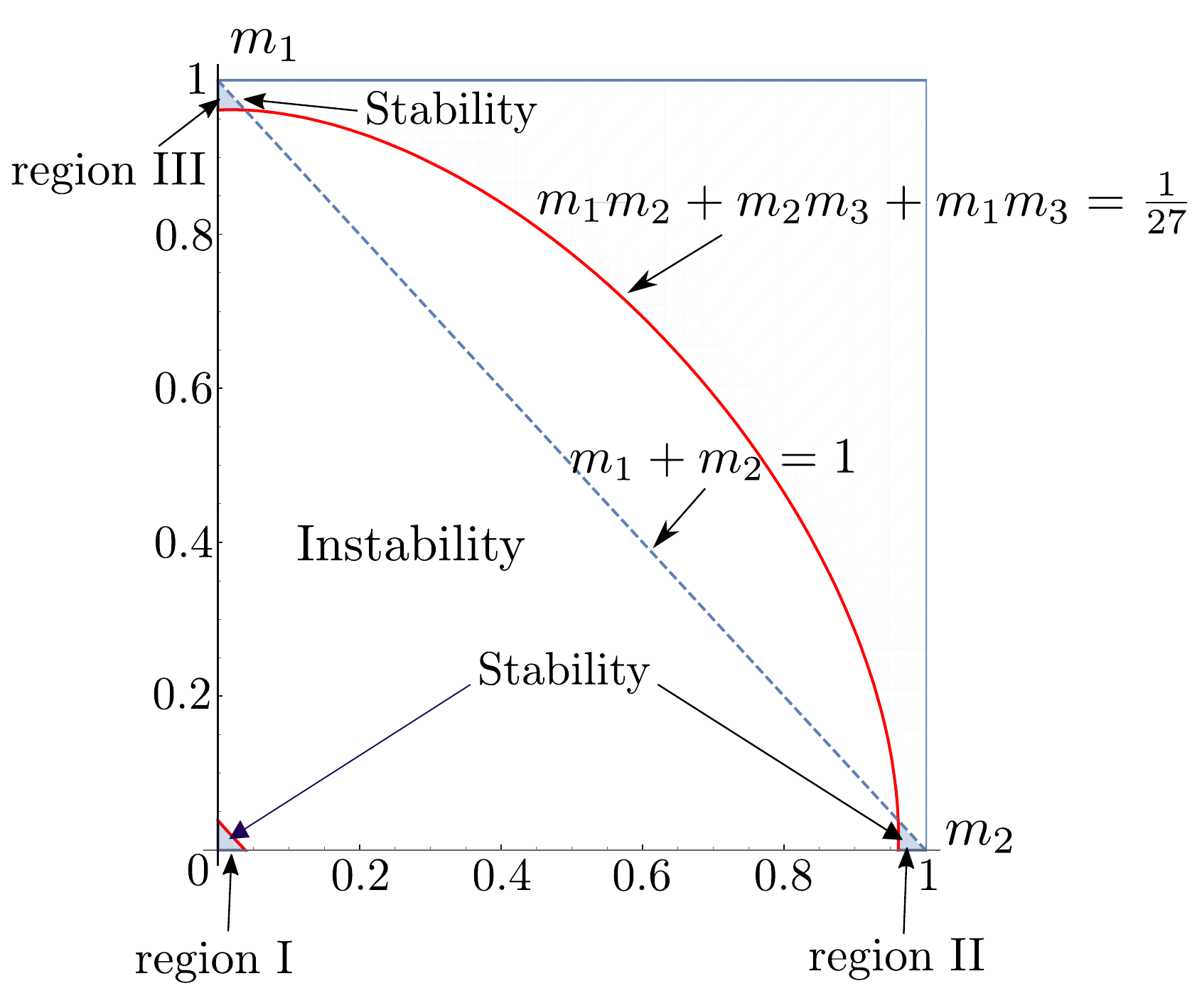}
\caption{The shaded triangular regions I, II, and III indicate the stability domains of
the Lagrange equilateral configuration of the primaries, while the white region
corresponds to instability. The third mass is given by
$m_3=1-m_1-m_2$, and the red curves represent Routh's critical curve.
Figure adapted from \cite{zepeAl3}.}\label{fig_routh}
\end{figure}

The total number of equilibrium points depends on the distribution of the
masses and may be eight, nine, or ten. The transition between the regions
with eight and ten equilibria occurs along a fold bifurcation curve, where
the system possesses exactly nine equilibrium points
\cite{pede, leandro, Balta1}. In the three linearly stable regions I,
II, and III of the admissible mass domain (Figure~\ref{fig_routh}), the
system has exactly eight equilibrium points.

In the present work, we restrict our attention to Regions I and III in the
$(m_1,m_2)$ plane (Figure~\ref{fig_routh}), where the Lagrange equilateral
configuration of the primaries is linearly stable. Because the three  regions are related by permutations of the mass parameters, they are
dynamically equivalent. We nevertheless present results for Regions I and III
to illustrate this symmetry and to provide a direct numerical comparison of the
corresponding stability diagrams.

\begin{figure}[H]
\centering
\includegraphics[scale=0.38]{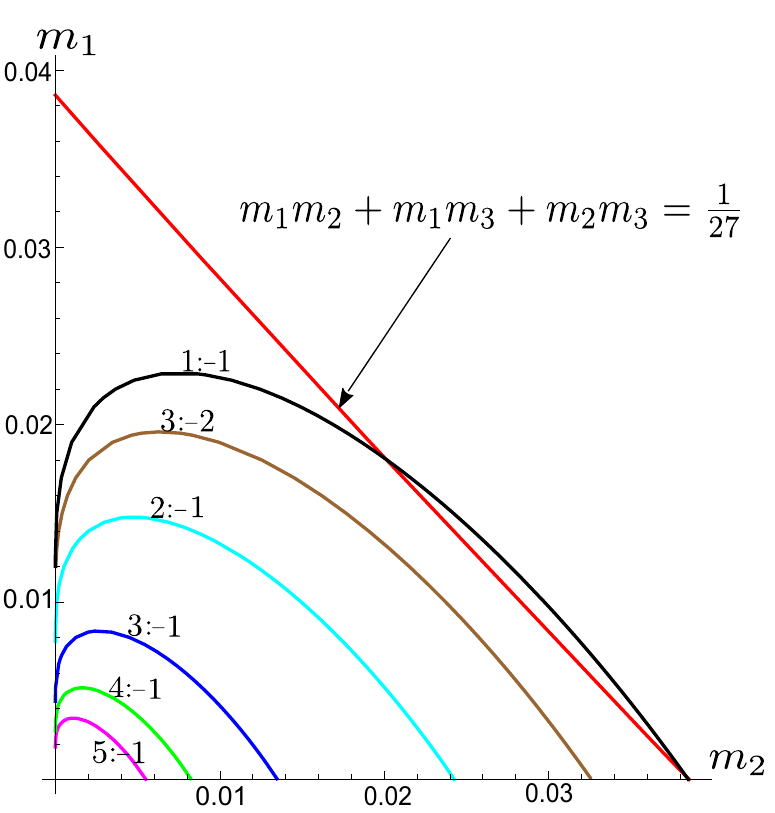}\qquad
\includegraphics[scale=0.38]{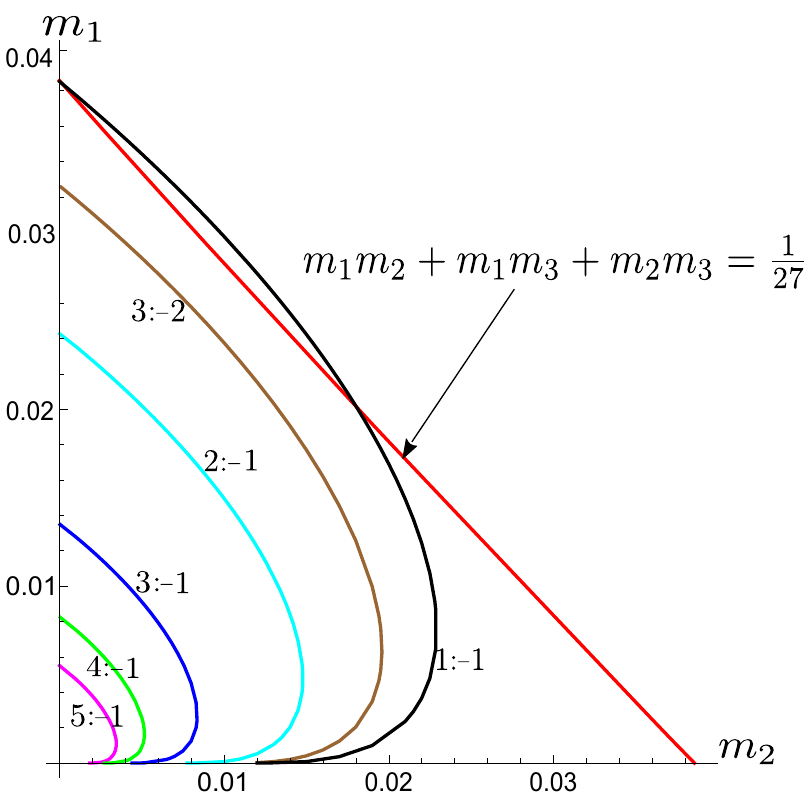}\qquad
\includegraphics[scale=0.38]{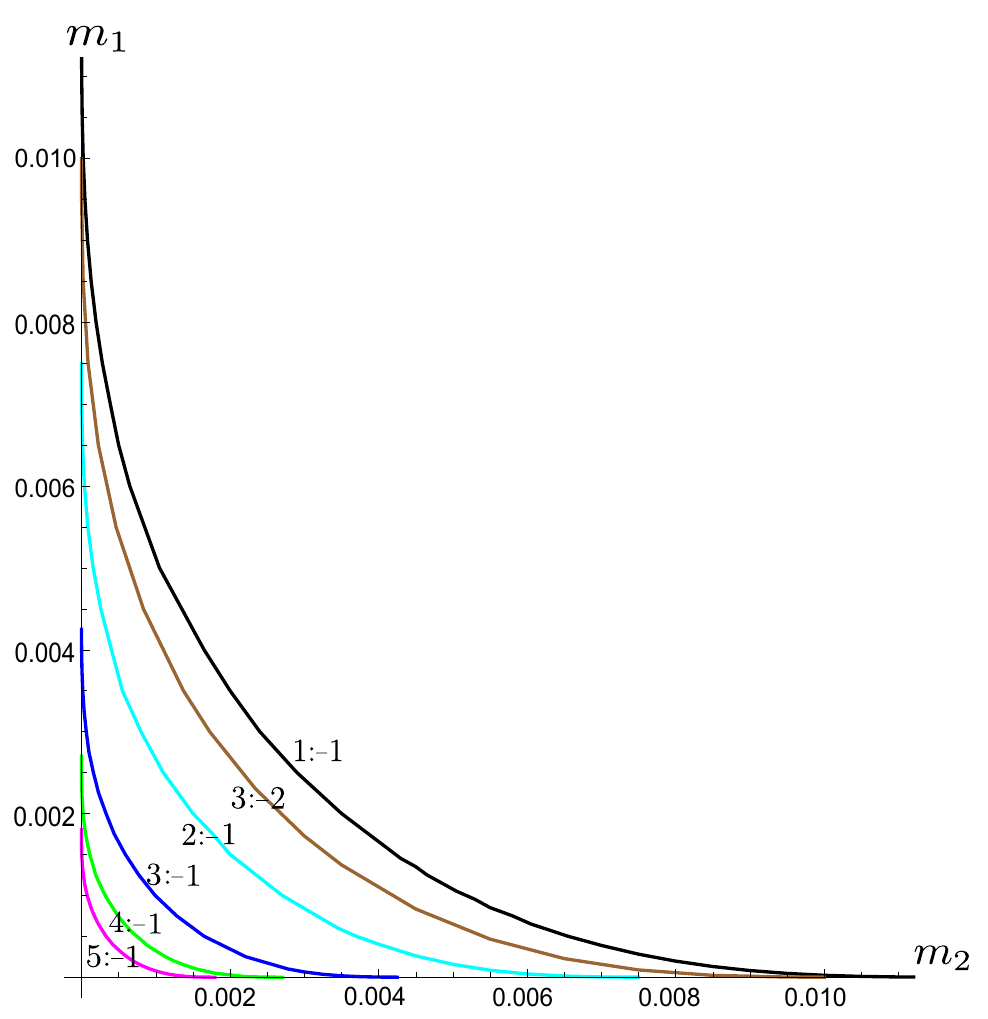}\qquad
\includegraphics[scale=0.38]{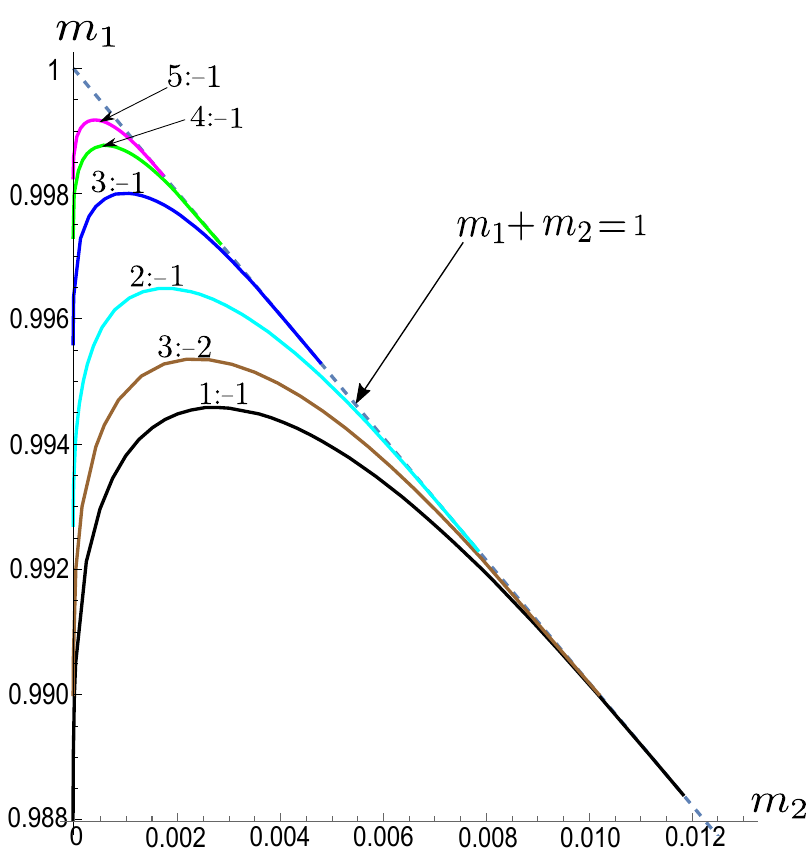}\qquad
\includegraphics[scale=0.38]{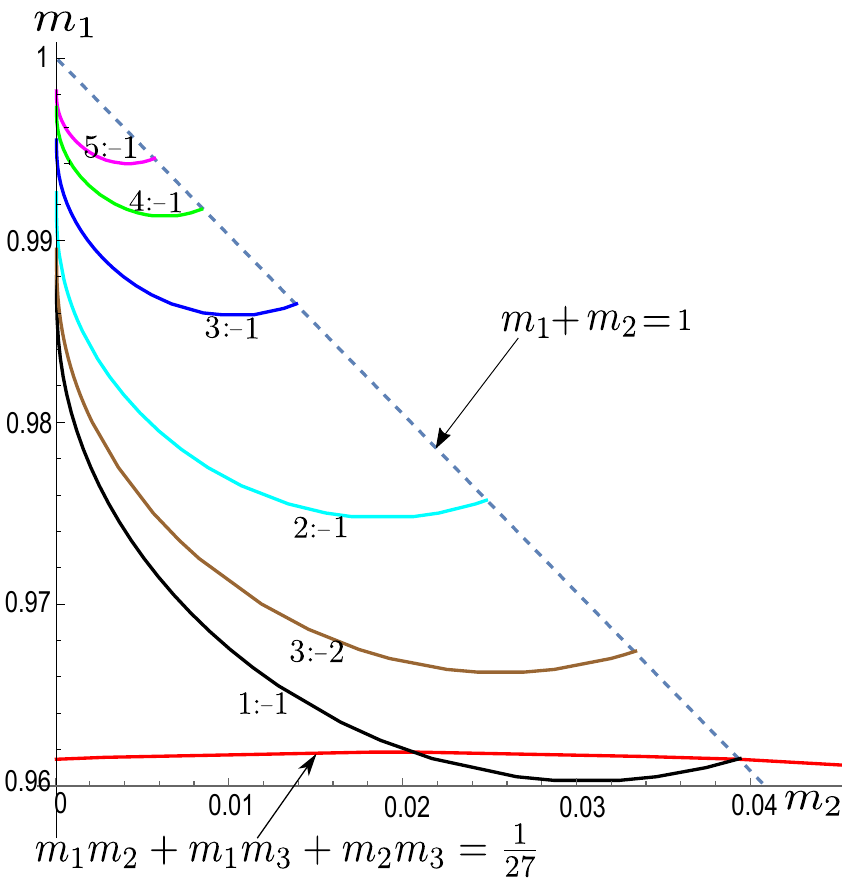}\qquad
\includegraphics[scale=0.38]{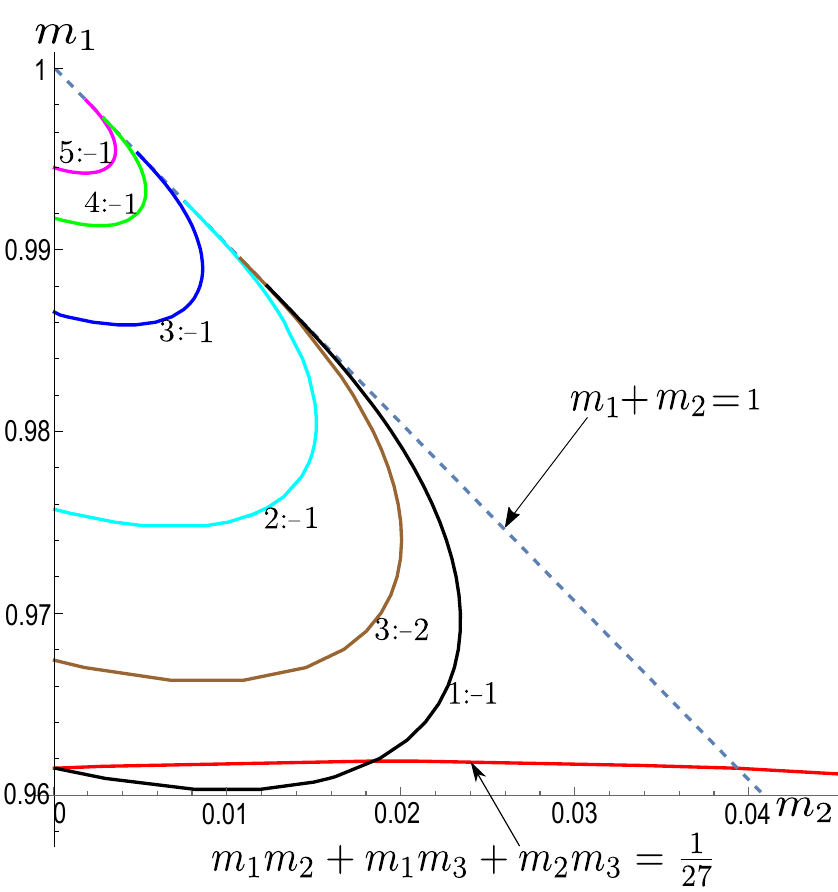}
\caption{Resonance curves associated with the elliptic equilibrium points
$L_3$, $L_5$, and $L_6$ in Regions I and III of the admissible
$(m_1,m_2)$ parameter plane introduced in Figure~\ref{fig_routh}. Adapted from \cite{zepeAl3}.}
\label{fig4}
\end{figure}

The linear stability of the equilibrium points has been extensively studied;
see, for instance,
\cite{pede,simo,zotos,Balta1,zepeAl2}.
These studies show that the points $L_1$, $L_2$, $L_4$, $L_7$, and $L_8$ are
linearly unstable, whereas $L_3$, $L_5$, and $L_6$ are elliptic in suitable
subregions of Regions I and III.
At the elliptic equilibria, the eigenvalues are of the form
\[
\lambda_{1,2}=\pm i\omega_1,
\qquad
\lambda_{3,4}=\pm i\omega_2.
\]
In \cite{zepeAl2}, resonance curves were constructed by imposing the condition
\[
k\omega_1(m_1,m_2)+m\omega_2(m_1,m_2)=0,
\]
where $k$ and $m$ are coprime integers.
For each elliptic equilibrium, several low-order resonance curves were obtained,
as shown in Figure~\ref{fig4}. The region of linear stability is bounded by the
$1\!:\!-1$ resonance curve, where $\omega_1=\omega_2$.

In the limiting case where one of the masses vanishes, the ERFBP reduces
to the classical restricted three-body problem. As shown in
Figure~\ref{fig4}, the resonance branches reach the boundary of the
admissible parameter domain at well-defined limiting values. Some of
these limiting values coincide with the classical critical mass ratios
of the restricted three-body problem, whereas others do not.
Consequently, the endpoint of a resonance branch in the ERFBP should not,
in general, be identified with a resonance of the limiting three-body
problem. This same phenomenon was already reported by Bardin and Volkov \cite{BardinVolkov2024} in
their numerical resonance diagrams.

Before proceeding to the nonlinear analysis, it is useful to summarize the
geometric information provided by Figures~\ref{fig_routh} and~\ref{fig4}, as
these figures establish the parameter regions and resonance structure used
throughout the remainder of the paper.

Figure~\ref{fig_routh} shows the mass plane $(m_1,m_2)$, with $m_3=1-m_1-m_2$
determined by the normalization~\eqref{eq_mass}. The red curve
$m_1m_2+m_1m_3+m_2m_3=\tfrac{1}{27}$ is Routh's critical curve: it separates
the shaded triangular regions I, II, and III, where the Lagrange equilateral
configuration of the primaries is linearly stable, from the white region,
where it is unstable. The dashed line $m_1+m_2=1$ (equivalently $m_3=0$)
marks the boundary of the admissible mass domain, along which the ERFBP
degenerates into the classical restricted three-body problem. 

Figure~\ref{fig4} restricts attention to Regions I and III and shows,
for each elliptic equilibrium $L_3$, $L_5$, and $L_6$, the loci in
$(m_1,m_2)$ where the resonance curves are displayed in each panel.
The outermost curve is always the $1$:$-1$ resonance, coinciding with
the boundary of the linear stability region; curves closer to the
origin correspond to higher-order resonances at smaller mass ratios.
Restricting attention to the $1\!:\!-1$ resonance curves in Region I, we find that
as either $m_1$ or $m_2$ tends to zero, one endpoint of two curves,
together with the endpoint of the remaining curve, approaches the
common limiting value
\[
m^*=
\frac12-\frac{\sqrt{15681+7104\sqrt{21}}}{450}
\approx0.01194203.
\]
This value is not associated with the classical $1$:$-1$ resonance
of the restricted three-body problem. Rather, it marks a fold
bifurcation of the ERFBP, at which four equilibrium points coalesce.
Indeed, as $m_2\to0$ with $m_1$ fixed and sufficiently small, the
triangular equilibrium $L_5$ of the restricted three-body problem
bifurcates into the four ERFBP equilibria $L_1$, $L_2$, $L_5$, and
$L_6$, which merge at $m^*$. Likewise, as $m_1\to0$ with $m_2$ small,
the triangular point $L_4$ gives rise to the four equilibria $L_3$,
$L_4$, $L_5$, and $L_7$, which also coalesce at the same parameter
value.
The opposite endpoints of the resonance branches associated with
$L_3$ and $L_6$ intersect the Routh stability boundary at
\[
m_R\approx0.03852089,
\]
which coincides with the classical Routh critical mass ratio of the
restricted three-body problem. In the corresponding limiting process,
$L_6$ converges to the triangular libration point $L_4$, whereas
$L_3$ converges to $L_5$. Consequently, these resonance branches
recover the classical $1$:$-1$ resonance of the restricted three-body
problem. This is consistent with the results reported in \cite{BardinVolkov2024}.
Therefore, Figure~\ref{fig4} identifies, prior to any normal-form computation,
which regions of the mass plane require the resonant analysis of
Section \ref{sec:birkhoffnorma}.

\section{Birkhoff normal form near elliptic equilibria}\label{sec:birkhoffnorma}
In order to analyze the nonlinear stability of the elliptic equilibrium points, we compute the normal form of the Hamiltonian in a neighborhood of each such point. In the case of elliptic equilibria, this normal form is known as the \emph{Birkhoff normal form}; see, for example, Theorem~5.5 in \cite{Arnold2006}.

To bring the system into normal form, we first shift the origin of the coordinate system to the equilibrium point $L_k = (x_0, y_0)$, where $k = 3, 5, 6$, by applying the linear change of variables defined by:$$\begin{array}{ll}\label{NL1}
	x=x_0+q_1, & \quad y=y_0+q_2,\\
	p_x=p_{x_0}+p_1, &  \quad p_y= p_{y_0}+p_2.
\end{array}$$
Next, we expand the Hamiltonian function \eqref{ham1} in a Taylor series around the origin. After dropping constant terms, we obtain a Hamiltonian of the form:
\begin{equation}\label{ham2}
	H=  H_2 + H_3+ H_4 + \cdots + H_m= \sum_{j=2}^m H_j,
\end{equation}
where each $H_j$ is a homogeneous polynomial of degree $ j$ in the variables $(q_1,q_2,p_1,p_2) $.

The second-order terms in the expansion of \eqref{ham2} are:
\begin{equation}
	H_2=\frac{1}{2}(p_1^2+p_2^2)+q_2p_1-q_1p_2+a_{20}q_1^2+a_{11}q_1q_2+a_{02}q_2^2,  \label{ham_cuad}
\end{equation}
with the associated Hamiltonian matrix
\begin{equation}\label{ham_matrix1}
Q=
\begin{pmatrix}
0 & 1 & 1 & 0\\
-1 & 0 & 0 & 1\\
2a_{20} & a_{11} & 0 & 1\\
a_{11} & 2a_{02} & -1 & 0
\end{pmatrix}.
\end{equation}

The four eigenvalues of the matrix $Q$ are $\pm i\omega_1$ and $\pm i\omega_2$, which are computed numerically using \emph{Mathematica}, considering mass values within the linear stability regions described in the previous section. 
We then compute the eigenvectors corresponding to the eigenvalues $+i\omega_1$ and $+i\omega_2$, denoted by $\alpha_1$ and $\alpha_2$, respectively. These eigenvectors are used to construct the coordinate transformation matrix ${\mathcal T}$, which maps the Cartesian coordinates $(q_1, q_2, p_1, p_2)$ to a new set of complex variables $(z_1, z_2, \bar{z}_1, \bar{z}_2)$. 
More specifically,
$$
{\mathcal T}=\left(\overline{\alpha_1}/r_1,\alpha_2/r_2,\alpha_1/r_1,\overline{\alpha_2}/r_2 \right),
$$
where $r_1$ and $r_2$ are constants necessary for ${\mathcal T}$ to be symplectic with symplectic multiplier $ \frac{1}{2}i$, so that ${\mathcal T}^TJ{\mathcal T}=\frac{1}{2}iJ$.  In terms of the complex coordinates, the quadratic part of the Hamiltonian \eqref{ham_cuad} takes the form
\begin{equation}
\overline{H}_2(z_1, \bar{z}_1, z_2, \bar{z}_2)=i\omega_1z_{1}\bar{z}_{1} -i\omega_2z_{2}\bar{z}_{2}\label{nfh0}, 
\end{equation}
and the full Hamiltonian \eqref{ham2} goes over
\begin{equation}\label{ham_com}
	{\overline H}(z_1,z_2,\bar{z}_1,\bar{z}_2) = {\overline H}_2 + 
	\sum_{\nu_1+\nu_2+\rho_1+\rho_2=3}^{m} h_{\nu_1\nu_2\rho_1\rho_2}z_1^{\nu_1}z_2^{\nu_2}\bar{z}_1^{\rho_1}\bar{z}_2^{\rho_2}.
\end{equation}

We note that the quadratic normal form ${\overline H}_2$ is not sign-definite; hence, Dirichlet’s theorem \cite{MeyerOffin} does not apply
 to determine the stability or instability of the equilibrium. Therefore, we proceed by normalizing the nonlinear part of \eqref{ham2} using the Lie-Deprit algorithm expressed in terms of the complex variables $(z_1, \bar{z}_1, z_2, \bar{z}_2)$. For details on its implementation, see \cite{zepAl} and \cite{MeyerOffin}.
This normalization procedure is applied exclusively to $H_3$, ensuring that each step of the algorithm produces a Hamiltonian that commutes with $\overline{H}_2$. The process is carried out up to a finite order; in most cases considered here, normalization up to fourth order is sufficient. This includes incorporating the quartic terms that define $H_4$.

To compute the normal form, we first determine the generating function associated with the symplectic transformation. Then, using \emph{Mathematica}, we numerically implement the Lie transformations to automate the normalization procedure. Due to the algebraic complexity and length of the resulting expressions, we omit the explicit form of the generating function.
Finally, we arrive at the normalized Hamiltonian in complex coordinates up to order four, which reads
\begin{equation}\label{com_realn}
	{\overline{\mathcal H}}(z_1,z_2,\bar{z}_1,\bar{z}_2)=i\omega_1z_{1}\bar{z}_{1}- i\omega_2z_{2}\bar{z}_{2}+i h_{2200} z_{1}^2\bar{z}_1^2 +
	i h_{1111} z_{1} z_2\bar{z}_1\bar{z}_2 + i h_{0022} z_{2}^2\bar{z}_2^2.
\end{equation}

In order to obtain the normalized Hamiltonian in real coordinates, we introduce the symplectic complex coordinate transformation defined as
\begin{equation}
	\begin{array}{lcl}
		z_{1}=\frac{1}{\sqrt{2}}(q_1+i p_1),  &\qquad &z_{2}=\frac{1}{\sqrt{2}}(q_2+i p_2),\vspace{0.3cm}\\
		\bar{z}_{1}=\frac{1}{\sqrt{2}}(q_1-i p_1),   &\qquad &
		\bar{z}_{2}=\frac{1}{\sqrt{2}}(q_2-i p_2),
	\end{array}\label{nf1b} 
\end{equation}
with  multiplier  $-\dfrac{1}{2}i$. Therefore, the normal form of the Hamiltonian \eqref{com_realn} in real coordinates becomes 
\begin{multline}
	\widetilde{\mathcal H}(q_1,q_2,p_1,p_2)=\frac{\omega_1}{2}\left( q^{2}_1+p^{2}_1\right) - \frac{\omega_2}{2}\left( q^{2}_2+p^{2}_2\right)+
	h_{2200}\left( q^{2}_1+p^{2}_1\right)^{2}\\
	+h_{1111}\left( q^{2}_1+p^{2}_1\right)\left(q^{2}_2+p^{2}_2\right)+h_{0022}\left( q^{2}_2+p^{2}_2\right)^2,\label{ham_realqp}
\end{multline}
where the coefficients 
 $h_{2200}$, $h_{1111}$   and  $h_{0022} $ are obtained numerically for each set of values of $ m_1 $, $ m_2 $, $ \omega_1 $ and $ \omega_2 $.

Assuming the presence of resonances, we proceed as before to construct the corresponding normal form. We begin by applying the Lie-Deprit algorithm to the Hamiltonian, normalized up to fourth-order terms and expressed in complex variables.

In the presence of a $2$:$-1$ resonance, that is, when $\omega_1 = 2\omega_2$, the normal form \eqref{ham2} reduces to
\begin{equation}
		{\overline{\mathcal H}} = i\left(2 z_{1}\bar{z}_{1}- z_{2}\bar{z}_{2}\right) 
		+ h_{1200}z_1z_2^{2} + h_{0012}\bar{z}_1\bar{z}_2^{2}.
		\label{nf_21}
\end{equation}
Next, we apply the transformation \eqref{nf1b} to express the Hamiltonian in Cartesian variables, which takes the form
\begin{multline}
	\widetilde{\mathcal H}(q_1,q_2,p_1,p_2)=\left( q^{2}_1+p^{2}_1\right) - \frac{1}{2}\left( q^{2}_2+p^{2}_2\right)+
	g_{1200} q_1q^{2}_2+g_{1002} q_1p^{2}_2+g_{0210}q_2^2p_1\\
	+g_{1101}q_1q_2p_2+g_{0111}q_2p_1p_2+g_{0012}p_1p_2^2. \label{ham_realqp1}
\end{multline}

Now we introduce the action-angle variables defined by  
\begin{equation}
I_j = \tfrac{1}{2}\,(q_j^2 + p_j^2), 
\qquad  
\phi_j = \tan^{-1}\!\left(\tfrac{p_j}{q_j}\right), 
\quad j = 1,2.
\end{equation}
Expressing the Hamiltonian ~\eqref{ham_realqp1} in these coordinates, we obtain
\begin{equation}\label{nf21}
{\mathscr H} 
= 2I_1 - I_2 
+ I_1^{1/2}I_2\big( A\cos (\phi_1+2\phi_2) + B\sin(\phi_1+2\phi_2)\big) 
= 2I_1 - I_2 + \delta I_1^{1/2}I_2 \cos\phi,
\end{equation}
where $\phi=\phi_1+2\phi_2$, $\delta=\sqrt{A^2+B^2}$, and $A$, $B$ are constants depending on the coefficients $h_{1200}$ and $h_{0012}$.

In contrast, under a $3$:$-1$ resonance, the normalized Hamiltonian expressed in complex variables takes the form:
\begin{equation*}\label{com_real21}
	{\overline{\mathcal H}} = i(3z_{1}\bar{z}_{1} -z_{2}\bar{z}_{2})+i h_{2200} z_{1}^2\bar{z}_1^2 +
	i h_{1111} z_{1} z_2\bar{z}_1\bar{z}_2 + i h_{0022} z_{2}^2\bar{z}_2^2+ h_{1300}z_1z_2^{3}+ h_{0013}\bar{z}_1\bar{z}_2^{3}.
\end{equation*}
Passing to Cartesian coordinates, this Hamiltonian becomes
\begin{multline}
	\widetilde{\mathcal H}(q_1,q_2,p_1,p_2)=\frac{3}{2}\left( q^{2}_1+p^{2}_1\right) - \frac{1}{2}\left( q^{2}_2+p^{2}_2\right)+h_{2200}\left( q^{2}_1+p^{2}_1\right)^{2}+h_{1111}\left( q^{2}_1+p^{2}_1\right)\left(q^{2}_2+p^{2}_2\right)
	\\
+h_{0022}\left( q^{2}_2+p^{2}_2\right)^2+g_{1300} q_1q^{3}_2+g_{1003} q_1p^{3}_2+g_{0310}q_2^3p_1+g_{1102}q_1q_2p_2^2\\
+g_{1201}q_1q_2^2p_2+g_{0112}q_2p_1p_2^2+g_{0211}q_2^2p_1p_2+g_{0013}p_1p_2^3,\label{ham_realqp2}
\end{multline}
which, in action-angle variables,  yields
\begin{align}\label{nf31}
	{\mathscr H} &=3I_1-I_2+ \frac{1}{2}(AI_1^{2}+BI_1I_2+CI^{2}_2)+I_1^{1/2}I_2^{3/2}(A_{13}\sin \phi_1+B_{13}\cos \phi_2)+\widetilde{H}\nonumber\\
	&=  3I_1-I_2+\frac{1}{2}(AI_1^{2}+BI_1I_2+CI^{2}_2)+I_1^{1/2}I_2^{3/2}\delta\cos \phi+\widetilde{H},
\end{align}
where
$\phi=\phi_1+3\phi_2$ and $\delta=\sqrt{A_{13}^2+B_{13}^2}$.
The constants \(A_{13}\) and \(B_{13}\) are determined by the coefficients
\(g\) of the resonant terms in the Hamiltonian \eqref{ham_realqp2}, and
\(\widetilde{H}\) denotes the collection of higher-order terms that are not
written explicitly.

\section{Nonlinear stability of the equilibria} \label{sec_nonstab}
This section presents a systematic numerical classification of the nonlinear
stability of the elliptic equilibria $L_3$, $L_5$, and $L_6$ throughout the
admissible mass plane. The analysis combines high-order Birkhoff normal forms
with the explicit computation of the Arnold determinants $D_4$ and $D_6$, together with
the numerical evaluation of the coefficients governing the low-order resonant
normal forms. Although the stability criteria employed below are classical,
their implementation in the unequal-mass equilateral restricted four-body
problem yields detailed stability boundaries, degeneracy loci, and resonance
structures that provide a global picture of the nonlinear dynamics.

In the absence of resonances between the frequencies $\omega_1$ and
$\omega_2$, the nonlinear stability of an elliptic equilibrium can be analyzed
using Arnold's stability theorem (Theorem~\ref{teo_arnold} in
Appendix~\ref{ap1}). To apply this result, we express the Birkhoff normal form
\eqref{ham_realqp} in action-angle variables:
\begin{equation}\label{ham_accion}
\mathscr{H}(I_1,I_2,\phi_1,\phi_2)
=
\omega_1 I_1-\omega_2 I_2
+
\frac{1}{2}
\left(
h_{2200}I_1^2
+
h_{1111}I_1I_2
+
h_{0022}I_2^2
\right).
\end{equation}

The nondegeneracy condition required by Arnold's theorem is expressed in terms
of the quartic invariant
\begin{equation}\label{eq_D4}
D_4=
h_{2200}\omega_1^2
+
h_{1111}\omega_1\omega_2
+
h_{0022}\omega_2^2.
\end{equation}
Whenever the frequencies are nonresonant and $D_4\neq 0$, Arnold's theorem
implies Lyapunov stability of the corresponding equilibrium point.

Since the equilibrium coordinates and the coefficients of the normal form are
not available in closed form, the nondegeneracy condition must be verified
numerically. For each equilibrium point, we construct a fine computational
grid in the admissible $(m_1,m_2)$ parameter domain, compute the associated
frequencies $\omega_1$ and $\omega_2$, detect low-order resonances, and
evaluate the corresponding invariants. This procedure produces detailed
numerical stability maps that identify the boundaries between stable,
resonant, and degenerate regimes.

As discussed in Section~\ref{sec:back}, the only equilibrium points that may
be elliptic are $L_3$, $L_5$, and $L_6$, located in Regions~I and~III of the
parameter plane.

\subsection{Region I}
\begin{figure}[h!]
\centering
\subfigure[]{
\includegraphics[width=0.50\textwidth]{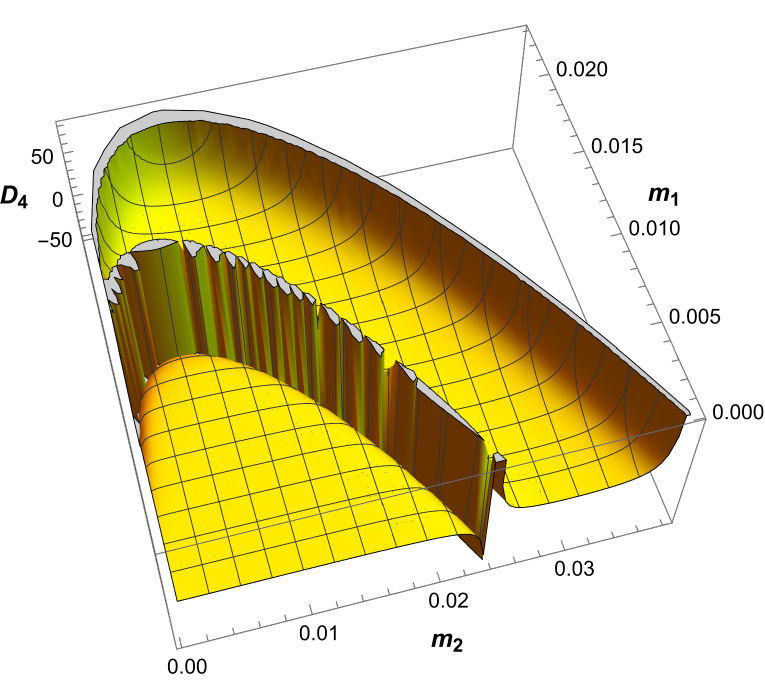}}
\hspace{0.6cm}
\subfigure[]{
\includegraphics[width=0.42\textwidth]{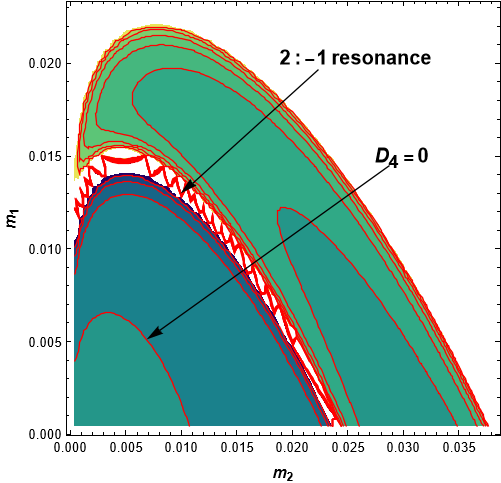}}
\caption{(a) Surface of $D_4$ for $L_3$ in Region~I.
(b) Corresponding contour plot. The curves $D_4=0$ and the
$2\!:\!-1$ resonance are indicated explicitly.}
\label{Enlf6}
\end{figure}

\begin{figure}[H]
\centering
\subfigure[]{
\includegraphics[width=0.45\textwidth]{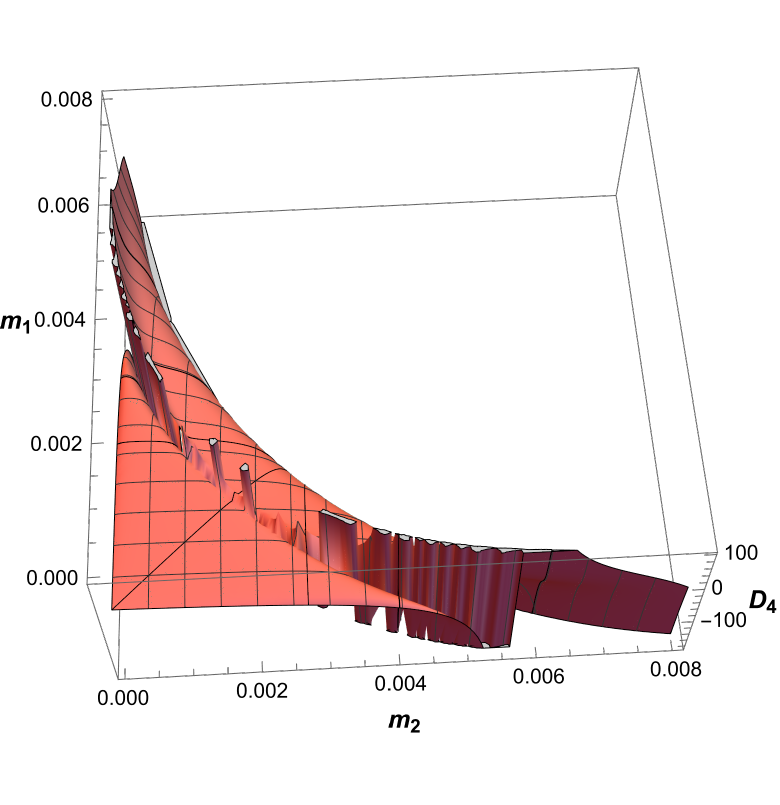}}
\hspace{0.6cm}
\subfigure[]{
\includegraphics[width=0.45\textwidth]{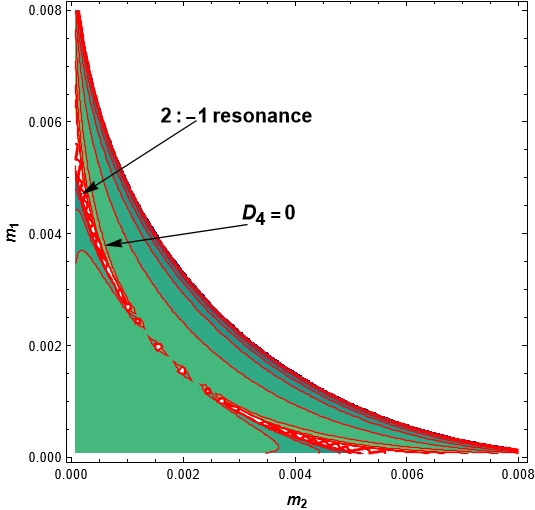}}
\caption{(a) Surface of $D_4$ for $L_5$ in Region~I.
(b) Corresponding contour plot.}
\label{Enlf8a}
\end{figure}

\begin{figure}[H]
\centering
\subfigure[]{
\includegraphics[width=0.45\textwidth]{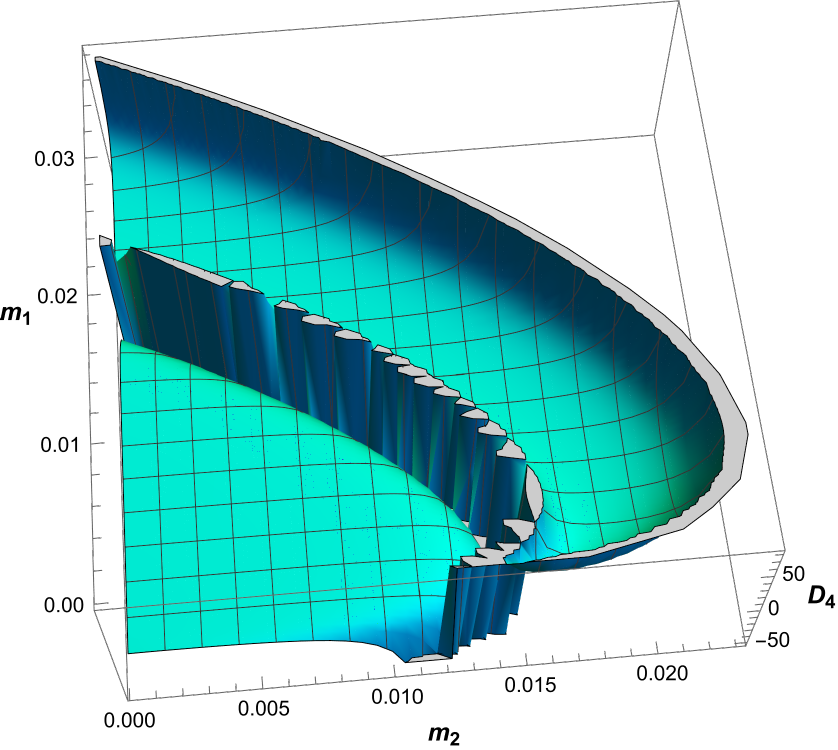}}
\hspace{0.4cm}
\subfigure[]{
\includegraphics[width=0.45\textwidth]{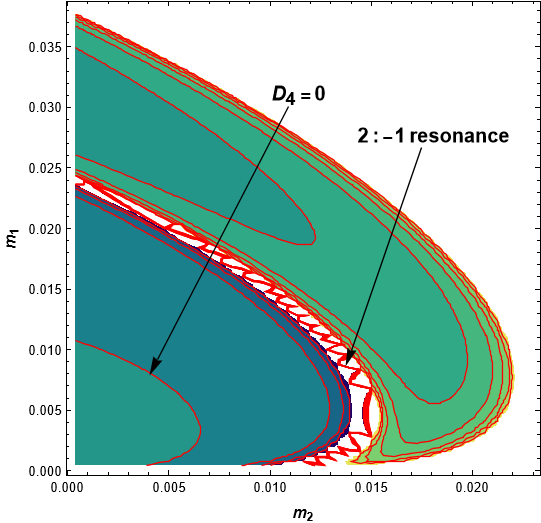}}
\caption{(a) Surface of $D_4$ for $L_6$ in Region~I.
(b) Corresponding contour plot.}
\label{Enlf8}
\end{figure}

Figures~\ref{Enlf6}(a), \ref{Enlf8a}(a), and \ref{Enlf8}(a) depict the
computed values of the quartic invariant $D_4$ for the equilibrium points
$L_3$, $L_5$, and $L_6$, respectively, throughout Region~I of the parameter
space. The corresponding contour plots are shown in
Figures~\ref{Enlf6}(b), \ref{Enlf8a}(b), and \ref{Enlf8}(b), where the curve
$D_4=0$ is explicitly indicated.

The discontinuities observed in these surfaces are numerically confirmed to
coincide with the $2\!:\!-1$ resonance curve, where the nonresonance
hypothesis required by Arnold's theorem fails. Away from the resonance and
degeneracy curves, the condition $D_4\neq 0$ holds, and Arnold's theorem
implies Lyapunov stability of the corresponding equilibrium point.

In particular, the critical configurations
$$(m_1,m_2)=(0.01,0.01495599), \quad \text{and} \quad (m_1,m_2)=(0.01,0.00006983),$$ 
previously identified in~\cite{zepAl}, lie on the curve
$D_4=0$ and are clearly visible in the contour plots.

\begin{figure}[H]
\centering
\subfigure[]{\includegraphics[scale=0.37]{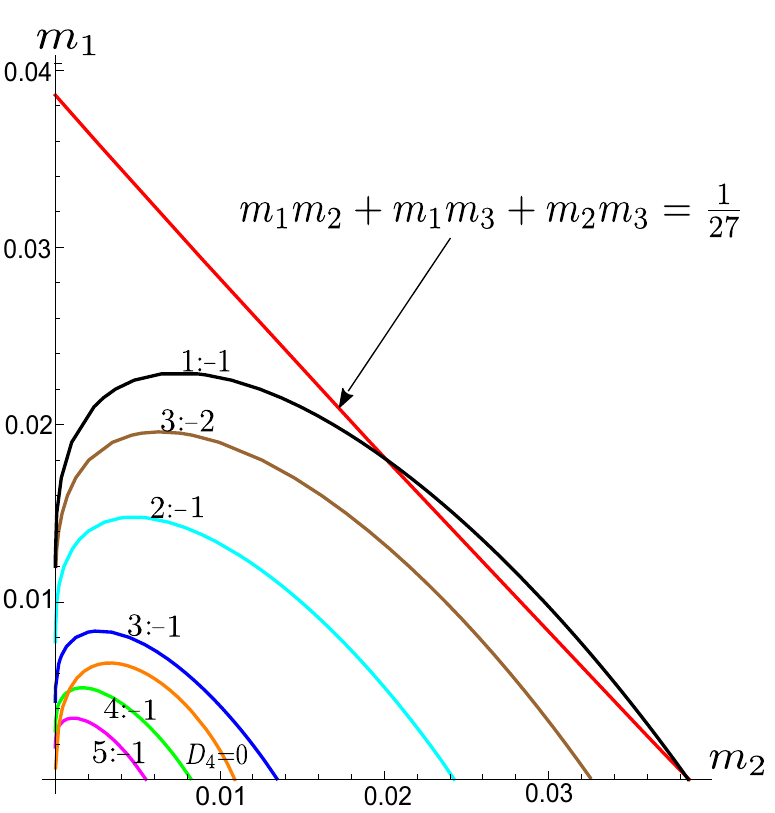}}
\hspace{0.8cm}
\subfigure[]{\includegraphics[scale=0.37]{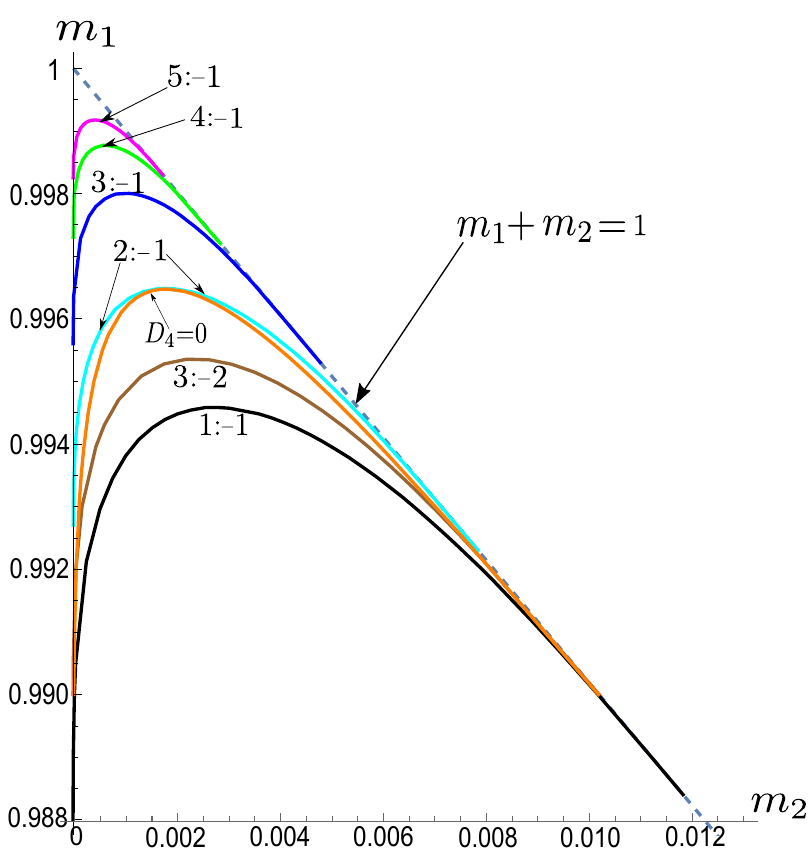}}
\hspace{0.8cm}
\subfigure[]{\includegraphics[scale=0.34]{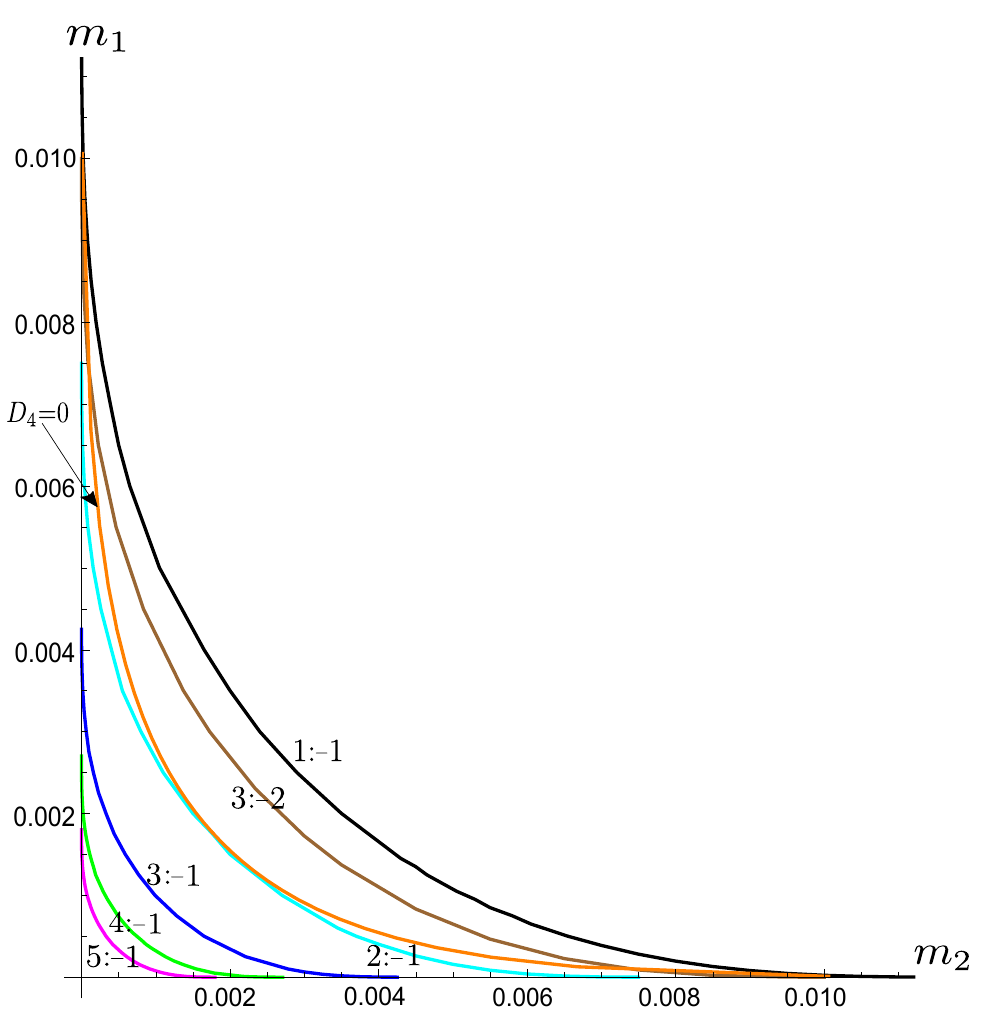}}
\hspace{0.8cm}
\subfigure[]{\includegraphics[scale=0.35]{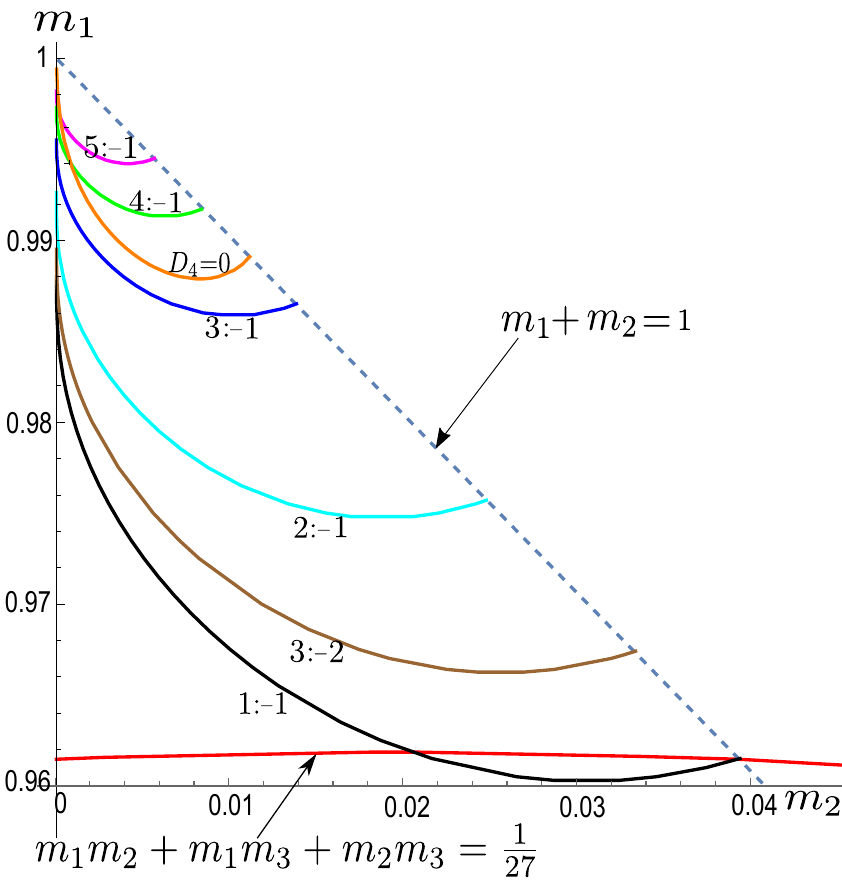}}
\subfigure[]{\includegraphics[scale=0.37]{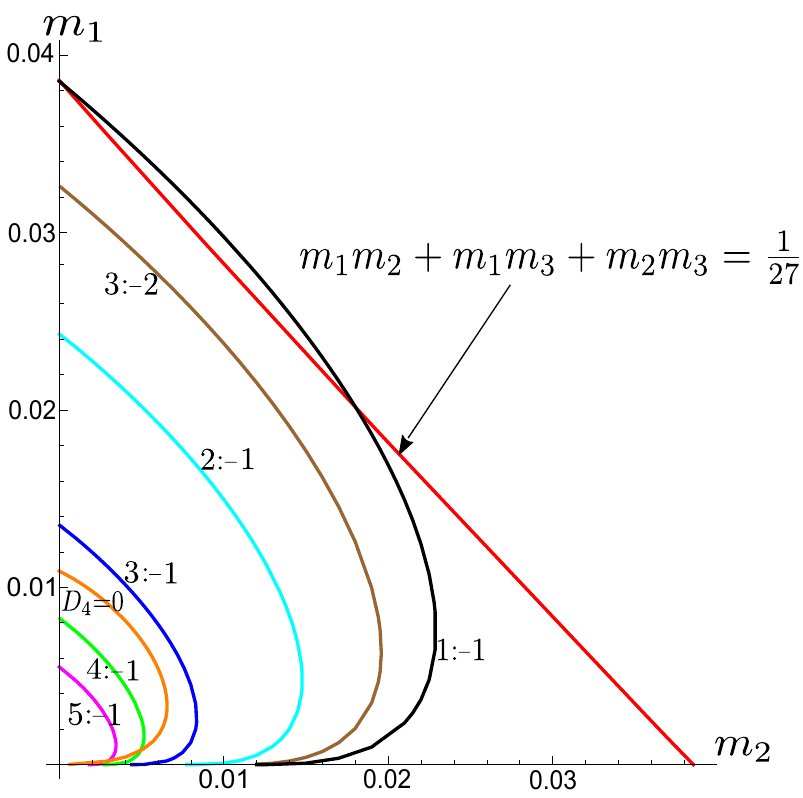}}
\hspace{0.8cm}
\subfigure[]{\includegraphics[scale=0.37]{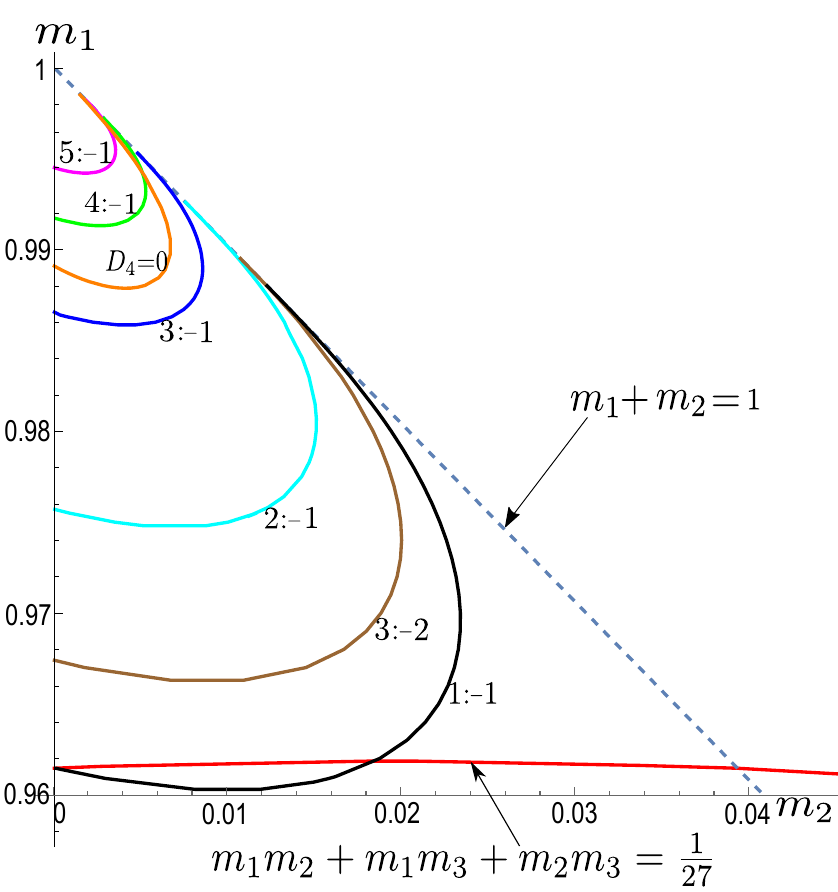}}
\caption{Resonance curves and the degeneracy curve $D_4=0$ (orange) for:
(a) $L_3$ in Region~I,
(b) $L_3$ in Region~III,
(c) $L_5$ in Region~I,
(d) $L_5$ in Region~III,
(e) $L_6$ in Region~I, and
(f) $L_6$ in Region~III.
Figures adapted from \cite{zepeAl3}.}
\label{d40res}
\end{figure}

\subsection{Region III}
A similar analysis is carried out for Region~III. The corresponding surfaces
and contour plots exhibit the same qualitative features: smooth variation of
$D_4$ in nonresonant domains, discontinuities along the $2\!:\!-1$
resonance curve, and well-defined curves where $D_4=0$.

Because Regions~I and~III are related by a permutation of the mass
parameters, they are dynamically equivalent. We nevertheless present
results for Regions I and III to illustrate this symmetry and to provide a
direct numerical comparison of the corresponding resonance structures.

\subsection{Stability at degenerate points ($D_4=0$)}
\label{sec:resona}
Figure~\ref{d40res} displays the numerically computed curves $D_4=0$
together with the relevant resonance curves for the equilibrium points
$L_3$, $L_5$, and $L_6$.

Along the curve $D_4=0$, the quartic nondegeneracy condition fails and the
fourth-order normal form is insufficient to determine stability. In this
case, the Lie--Deprit normalization is extended to sixth order:
\begin{multline*}
\mathscr{H}
=
\omega_1 I_1-\omega_2 I_2
+
\frac{1}{2}
(h_{2200}I_1^2+h_{1111}I_1I_2+h_{0022}I_2^2)
\\
+
\frac{1}{2}
(h_{3300}I_1^3+h_{2211}I_1^2I_2
+h_{1122}I_1I_2^2+h_{0033}I_2^3).
\end{multline*}

The corresponding sixth-order invariant is
\[
D_6
=
h_{3300}\omega_2^3
+
h_{2211}\omega_2^2\omega_1
+
h_{1122}\omega_2\omega_1^2
+
h_{0033}\omega_1^3.
\]

The invariant $D_6$ is evaluated numerically along the degeneracy curve.
At parameter values where $D_6\neq 0$, the sixth-order terms provide the
leading nontrivial contribution to the normal form and yield the
corresponding nonlinear stability criterion. Parameter values for which
$D_6=0$ correspond to higher-order degeneracies and are not considered in
the present work.

Similar higher-order degeneracies are known to occur in the classical
restricted three-body problem; see, for example,
\cite{Markeev,MeyerOffin}. The resulting stability classification is summarized in
Table~\ref{res_tabr1}.
\begin{table}[h!]
\centering
\begin{tabular}{|c|c|c|}
\hline
Equilibrium point & Region & Stability\\
\hline
$L_3$ & I   & stable \\ \hline
$L_5$ & I   & stable \\ \hline
$L_6$ & I   & stable \\ \hline
$L_3$ & III & stable \\ \hline
$L_5$ & III & stable \\ \hline
$L_6$ & III & stable \\ \hline
\end{tabular}
\caption{Stability of the elliptic equilibria along the degeneracy curve
$D_4=0$ at parameter values where $D_6\neq 0$ and no resonance is present.}
\label{res_tabr1}
\end{table}
Alvarez-Ramírez et al.~\cite{alvastuchi} previously analyzed the special case
of two equal masses. Their results are recovered as particular cases of the
general framework developed here and are summarized in Table~\ref{tab_vier}.
\begin{table}[H]
\centering
\begin{tabular}{|c|c|c|c|c|}
\hline
Region & Equal masses & Equilibrium point & $D_6$ & Stability \\
\hline
I  & 0.005809004 & $L_3$, $L_6$ & 546.371       & stable \\ \hline
I  & 0.00176397  & $L_5$        & 1574120.9148 & stable \\ \hline
III & 0.00176397 & $L_3$        & 1574120.9148 & stable \\ \hline
III & 0.005809004 & $L_5$, $L_6$ & 546.371      & stable \\ \hline
\end{tabular}
\caption{Stability of $L_3$, $L_5$, and $L_6$ for configurations with two
equal masses, where $D_4=0$ and $D_6\neq 0$.}
\label{tab_vier}
\end{table}

\subsection{Nonlinear behavior under the $2\!:\!-1$ resonance}
\label{sec:resona21}
Since Arnold's theorem does not apply in the resonant case, we use the Alfriend-Markeev-Meyer theorem (Theorem~\ref{teo_AMM} in
Appendix~\ref{ap1}).

The coefficient $\delta$ appearing in the resonant normal form
\eqref{nf21} is computed numerically using variable-precision arithmetic.
Its value was verified at representative points along the resonance curve,
and in all cases considered we found $\delta\neq 0$.
The resulting stability classification is summarized in Table~\ref{res_tabr2}.
\begin{table}[hbtp]
\centering
\begin{tabular}{|c|c|c|}
\hline
Equilibrium point & Region & Stability\\
\hline
$L_3$ & I   & unstable \\ \hline
$L_5$ & I   & unstable \\ \hline
$L_6$ & I   & unstable \\ \hline
$L_3$ & III & unstable \\ \hline
$L_5$ & III & unstable \\ \hline
$L_6$ & III & unstable \\ \hline
\end{tabular}
\caption{Stability of the elliptic equilibria along the $2\!:\!-1$
resonance curve.}
\label{res_tabr2}
\end{table}

Table~\ref{tablemm} lists the corresponding results for configurations with
two equal masses.

\begin{table}[hbt!]
\centering
\begin{tabular}{|c|c|c|c|c|}
\hline
Region & Equal masses & Equilibrium & $\delta$ & Stability \\
\hline
I   & 0.01208   & $L_3$, $L_6$ & 4.19339  & unstable \\ \hline
I   & 0.00100005 & $L_5$       & 0.725005 & unstable \\ \hline
III & 0.00679   & $L_3$       & 0.725005 & unstable \\ \hline
III & 0.01208   & $L_5$, $L_6$ & 4.19339 & unstable \\ \hline
\end{tabular}
\caption{Stability of $L_3$, $L_5$, and $L_6$ under the $2\!:\!-1$
resonance for equal-mass configurations.}
\label{tablemm}
\end{table}

\subsection{Nonlinear stability in the case of $3\!:\!-1$ resonance}
\label{sec:resona31}
To analyze the $3\!:\!-1$ resonance, we apply Theorem~\ref{teo31} to the
normal form \eqref{nf31}. The relevant quantity is
\[
\Delta D = 6\sqrt{3}\,|\delta| - |D|,
\qquad
D = A + 6B + 9C.
\]

According to Theorem~\ref{teo31}, the equilibrium is stable if
$\Delta D<0$, unstable if $\Delta D>0$, and the criterion is inconclusive
when $\Delta D=0$.

The numerical evaluation of $\Delta D$ was carried out using
\textit{Mathematica} along the resonance curve.

In Region~I, the equilibria $L_3$ and $L_6$ are unstable along the entire
$3\!:\!-1$ resonance curve, whereas $L_5$ is stable on a finite segment and
unstable outside it. In Region~III, the equilibria $L_3$ and $L_6$ are again
unstable along the whole resonance curve, while $L_5$ exhibits the
complementary behavior.

Representative values for equal-mass configurations are listed in
Table~\ref{tab_nn}.

\begin{table}[hbt!]
\centering
\begin{tabular}{|c|c|c|c|c|}
\hline
Region & Equal masses & Equilibrium Point & $\Delta D$ & Stability \\
\hline
I   & 0.00679    & $L_3$, $L_6$ & 440.763     & unstable \\ \hline
I   & 0.00100005 & $L_5$        & -8276.4446  & stable \\ \hline
III & 0.00100005 & $L_3$        & -8276.4446  & stable \\ \hline
III & 0.00679    & $L_5$, $L_6$ & 440.763     & unstable \\ \hline
\end{tabular}
\caption{Stability of $L_3$, $L_5$, and $L_6$ under the $3\!:\!-1$
resonance for equal-mass configurations.}
\label{tab_nn}
\end{table}
To the best of our knowledge, this appears to be the first explicit numerical
and analytical classification of the nonlinear behavior associated with the
$3\!:\!-1$ resonance in the unequal-mass ERFBP, including the resulting
stability and instability properties of the equilibrium points
$L_3$, $L_5$, and $L_6$.

The numerical results obtained in this section show that the nonlinear
dynamics of the elliptic equilibria are organized by three fundamental
structures:
\begin{enumerate}
\item[{1.}] Nonresonant regions where $D_4\neq 0$ and Arnold's theorem guarantees
Lyapunov stability.

\item[{2.}] Degeneracy curves where $D_4=0$ and sixth-order terms, through the
invariant $D_6$, determine the leading stability properties.

\item[{3.}] Resonance curves, particularly the $2\!:\!-1$ and $3\!:\!-1$
commensurabilities, where resonant normal forms govern the local dynamics.
\end{enumerate}
Taken together, these computations provide a systematic numerical
classification of the nonlinear stability of the elliptic equilibria
$L_3$, $L_5$, and $L_6$. In Section~\ref{sec:singular}, these results are
reinterpreted geometrically through singular reduction, which reveals the
bifurcation mechanisms underlying the observed transitions between stability
and instability.

\section{Singular reduction and bifurcation geometry}\label{sec:singular}
The analysis of Section \ref{sec_nonstab} established the nonlinear stability of $L_3$, $L_5$, and $L_6$
  in both nonresonant and resonant regimes using Birkhoff 
normal form and the theorems of Arnold, Alfriend-Markeev-Meyer and  Alfriend-Markeev. However, near 
the $k\!:\!-1$ resonance manifolds, this approach remains essentially local, since 
the dynamics is further structured by invariants associated with rotational symmetry. 
To capture the global geometry of resonant dynamics, we now turn to singular 
reduction. In \cite{palacian2018}, Meyer et al.\ provide a geometric interpretation of singular
reduction theory in the context of resonant Hamiltonians with $n$ degrees of freedom. 
In addition, Section~12.6 of \cite{MeyerOffin} illustrates this technique in the restricted three-body problem. 
Here, we extend these ideas to the ERFBP, reducing with respect to the rotation generated by the resonant angle 
$\phi_1 + k\phi_2$ in order to analyze the existence, stability, and bifurcations of periodic orbits.

Let us consider the quadratic part of the Hamiltonian, expressed in normal form,  under resonance
\begin{equation}\label{cuadratico}
\mathbb{H}=\frac{1}{2}\left[ k(q_1^2+p_1^2)-(q_2^2+p_2^2)\right]=kI_1-I_2,
\end{equation}
where $k  =\dfrac{ \omega_1}{\omega_2}= 2, 3$, $r = (q_1, q_2, p_1, p_2)$ is the position vector in Cartesian coordinates, and $I_1, I_2, \phi_1, \phi_2$ are action-angle variables. The associated Hamiltonian system is given by:
\begin{equation}\label{cap112}
	\dot{I}_1 = 0, \qquad \dot{I}_2 = 0, \qquad \dot{\phi}_1 = -k, \qquad \dot{\phi}_2 = 1.
\end{equation}
This system admits three  independent integrals of motion:  $I_1$, $I_2$, and the linear combination $\phi_1 + k\phi_2$.
These invariants are sufficient to fully characterize the orbit dynamics in the four-dimensional phase space.

The fundamental set of invariant polynomials associated with the $k$:$-1$ resonance is
\begin{equation} \label{rs1}
	\begin{split}
		a_1& =I_1=q_1^2+p_1^2,   \\
		a_2&=I_2=q_2^2+p_2^2,  \\
		a_3&=I_1^{1/2}I_2^{k/2}\cos(\phi_1+k\phi_2)=\mbox{Re}[(q_1+ip_1)(q_2+ip_2)^k],  \\
		a_4&=I_1^{1/2}I_2^{k/2}\sin(\phi_1+k\phi_2)=\mbox{Im}[(q_1+i p_1)(q_2+ip_2)^k].
	\end{split}
\end{equation}
These invariants are not independent but satisfy the algebraic relation
\begin{equation} \label{eq:relation}
a^2_3 + a^2_4 = a_1 a^k_2,
\end{equation}
which is a direct consequence of the Pythagorean identity 
$\cos^{2}\phi + \sin^{2}\phi = 1$ 
together with the nonnegativity conditions $a_{1} \geq 0$ and $a_{2} \geq 0$.

The Poisson brackets among these invariants are listed in Table~\ref{tabrs1}. 
Because $k$ is a positive integer (in particular $k=2,3$ in the cases of interest), 
all brackets take polynomial form. After presenting the general structure of the 
invariants and their Poisson brackets for arbitrary $k$, we now focus on the 
resonances $k=2$ and $k=3$. These constitute the simplest nontrivial cases: the 
case $k=2$ admits an explicit analysis of the reduced dynamics, while $k=3$ 
illustrates the next level of complexity and the extension of the algebraic 
structure. Other values of $k$ can be handled in a similar way.
\begin{table}[hbt]
	\centering 
	\begin{tabular}{| c | c | c | c |c|}
		\hline 
		$ \{ \cdot , \cdot \}  $	& $a_{1}  $ & $ a_{2} $ & $ a_{3} $&$ a_{4} $ \\ 
		\hline
		$ a_{1} $	& 0 &0  &$- 2a_{4} $ &  $ 2a_{3} $\\ \hline
		$ a_{2} $	&0  &0  &$- 2ka_{4} $&$ 2ka_{3} $\\ \hline
		$ a_{3} $	& $ 2a_{4} $ &$ 2ka_{4} $  &0 &$ a_{2}^{k-1}(k^{2}a_{1}+a_{2}) $ \\ \hline
		$ a_{4} $	& $- 2a_{3} $ &$- 2ka_{3} $ &$- a_{2}^{k-1}(k^{2}a_{1}+a_{2}) $& 0\\  
		\hline 
	\end{tabular}
	\caption{Poisson brackets of the invariants $a_1$, $a_2$, $a_3$, $a_4$.}\label{tabrs1}
\end{table}

Since the Hamiltonians \eqref{nf21} and \eqref{nf31} are already in normal form, the quantity   $\mathbb{H}$ is conserved.  
We therefore fix its value by setting
$$h = \mathbb{H} = k I_1 - I_2 = k a_1 - a_2,$$
corresponding to $k=2,3$, respectively.  
From this relation one obtains the dependent variables as
$$
a_2 = k a_1 - h,  \qquad  I_2 = k I_1 - h,
$$
and analogous expressions follow.
\medskip

\noindent\textbf{Case} $\mathbf{k=2}$: 
In the $2$:$-1$ resonance, the Hamiltonian \eqref{nf21} can be written in terms of the invariant polynomials as
$${\mathscr H} =2a_1 - a_2 + \delta a_3 = h + \delta a_3.$$
Discarding the constant term $h$, the reduced averaged Hamiltonian simplifies to  $\bar{R} = \delta a_3$.

In Section~\ref{sec_nonstab}, it was shown numerically that the set of $(m_1,m_2)$ values for which the frequencies 
$\omega_1$ and $\omega_2$ of the ERFBP linearization around the points $L_3$, $L_5$, and $L_6$ 
satisfy the $2\!:\!-1$ resonance condition forms a continuous curve. 
Our numerical computations further establish that $\delta \neq 0$  for all admissible 
mass values along this curve. After an appropriate rescaling of time, we may 
normalize to $\delta=1$, so that the reduced averaged Hamiltonian simplifies to $\bar{R} = a_3.$
Consequently, singular reduction theory applies uniformly along the entire resonance curve.

The geometric structure of the problem is determined by the intersection of two surfaces in $\mathbb{R}^3$. 
The first one, the \emph{orbit space} $\mathbb{O}$, is defined by the invariant relation
$$
a^2_3+a_4^2=a_1(2a_1-h)^2,\quad a_1\geq0, \quad a_1\geq\frac{1}{2}h,
$$
The second surface is associated with the reduced averaged Hamiltonian and is given by
$$
\bar{R}=\bar{h}.
$$
According to the value of $\bar{h}$, the orbit space near the origin in the $(a_1,a_3,a_4)$ coordinate system assumes different shapes:  
a paraboloid of revolution when $\bar{h}<0$, yielding a smooth surface;   a rotated cusp when $\bar{h}=0$;  
and a cone when $\bar{h}>0$.

\begin{figure}[!h]
	\subfigure[]{\includegraphics[width=0.32\textwidth]{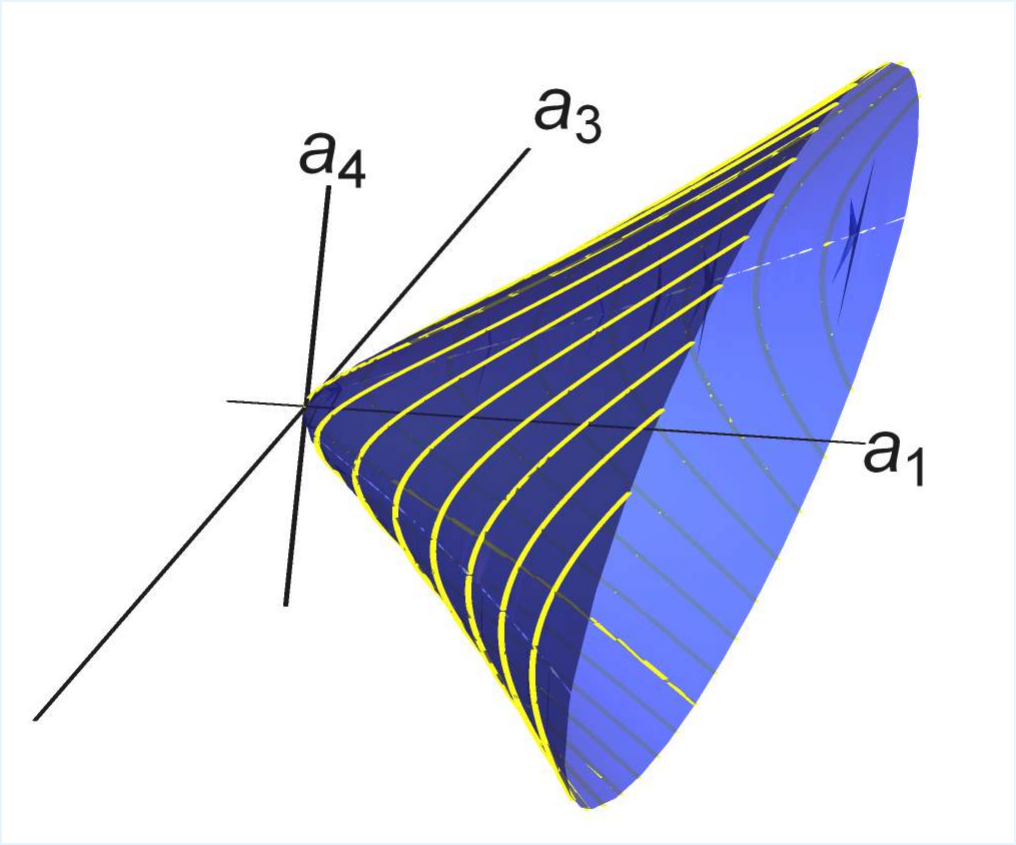}}
	\subfigure[]{\includegraphics[width=0.32\textwidth]{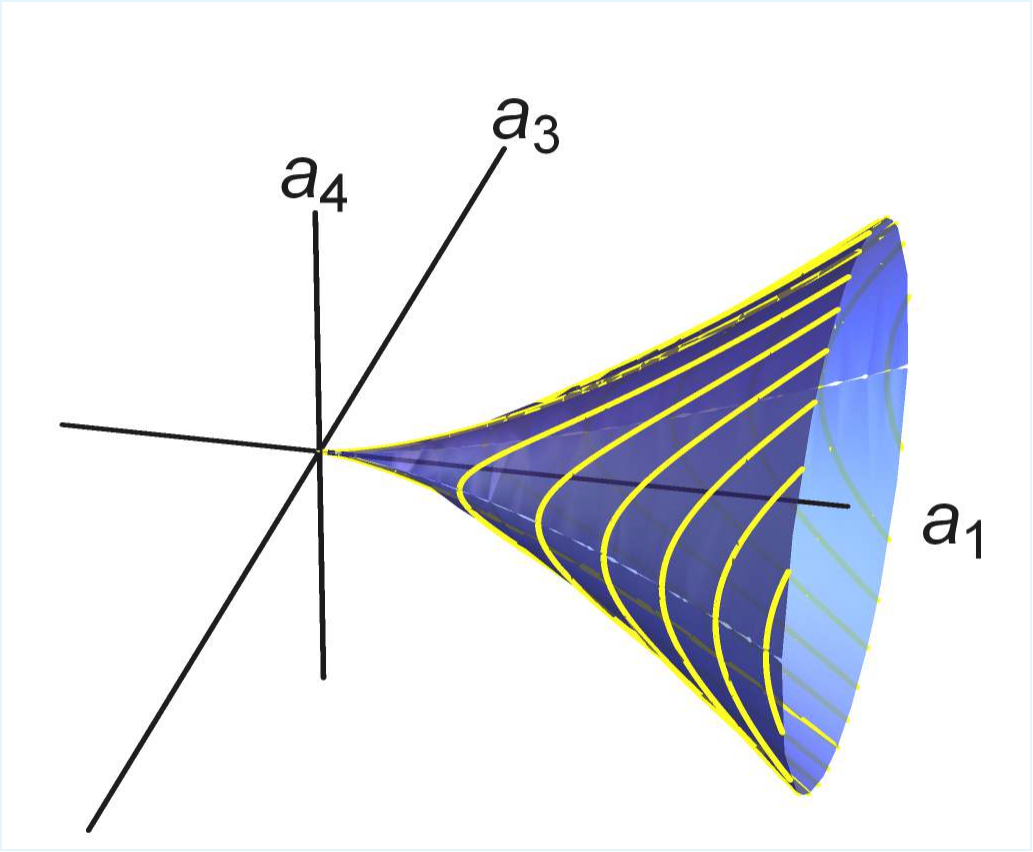}}
	\subfigure[]{\includegraphics[width=0.32\textwidth]{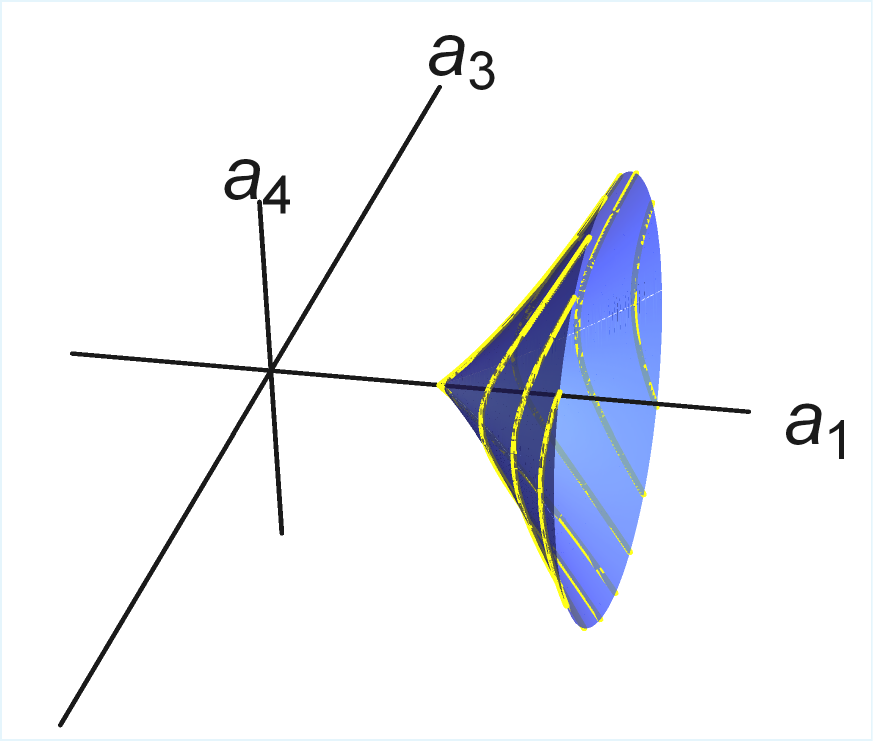}}	
	\caption{Flow on the orbit space ${\mathbb O}$ for the $2$:$-1$ resonance where (a) $h<0$ , (b) $h=0$ and (c) $h>0$.} \label{rs3a}
\end{figure}

From a geometric viewpoint, the constraint $\bar{R} = a_3 = \text{const}$ defines a  parallel plane to the $(a_1,a_4)$-plane.  
The induced vector field determines the dynamics, that is, the flow restricted to the orbit space is described by the corresponding equations of motion:
$$
\dot{a}_1=\{a_1,\bar{R} \}=-2a_4,\quad\dot{a}_3=\{a_3,\bar{R}\}=0,\quad\dot{a}_4=\{a_4,\bar{R}\}=-a_2(4a_1+a_2).
$$
Figure~\ref{rs3a} illustrates the flow on the orbit space ${\mathbb O}$ for the cases $h<0$, $h=0$, and $h>0$.

The equilibrium solutions follow from the conditions $a_4=0$ and $a_2(4a_1+a_2)=0$.  
Because $a_1$ and $a_2$ cannot simultaneously vanish or take negative values, owing to the relation $a_2 = 2a_1 - h$,  
the only admissible equilibrium configuration is given by $a_4 = a_2 = 0$.  
This requires $a_2=0$, which is satisfied when $a_1 = h/2$, and consequently implies $h \geq 0$.
  
The dynamics on the orbit spaces, illustrated in Figure~\ref{rs3a}, takes place on the invariant planes $\bar{R} = a_3 = \text{const}$, with the coordinate $a_4$ strictly decreasing along the flow since $\dot{a}_4 < 0$.  
The qualitative behavior of the system depends on the sign of the parameter $h$, leading to three distinct geometric scenarios.  
When $h<0$, the orbit space contains no equilibria. As a result, all trajectories escape as $t \to \pm \infty$, and no periodic families are generated.  
At the degenerate case $h=0$, the orbit space passes through the origin of $\mathbb{R}^3$. The dynamics exhibits two heteroclinic-like trajectories in $\mathbb{O}$ approaching the origin in forward and backward time, indicating that the equilibrium is unstable. 
For $h>0$, a bifurcation occurs: an equilibrium appears at $(a_1,a_3,a_4) = (h/2,0,0)$. This point corresponds to a Lyapunov family of short-period orbits with period $T \sim \pi$, for each $h \geq 0$. Furthermore, the existence of an orbit asymptotically approaching this equilibrium as $t \to \pm \infty$ reveals that this family is unstable.  
In summary, the transition across $h=0$ separates three regimes: divergence without equilibria ($h<0$), a degenerate unstable configuration ($h=0$), and the appearance of the short-period Lyapunov family ($h>0$).

\medskip

\noindent\textbf{Case} $\mathbf{k=3}$.  
The Hamiltonian \eqref{nf31} expanded up to order four is
$${\mathscr H} = 3I_1 - I_2 + \frac{1}{2}(A I_1^{2} + 2 B I_1 I_2 + C I_2^{2})
+ \delta I_1^{1/2} I_2^{3/2} \cos \phi,$$
where $A$, $B$, $C$, and $\delta$ are constants determined by the fundamental frequencies $\omega_{1}(m_1,m_2)$ and $\omega_{2}(m_1,m_2)$.  
In terms of the invariants, this expression becomes
\begin{equation}\label{hred31}
{\mathscr H} = 3a_1 - a_2 + \frac{1}{2}(A a_1^2 + 2B a_1 a_2 + C a_2^2) + \delta a_3.
\end{equation}

Fixing the energy level in the Hamiltonian $\mathbb{H} = 3 I_1 - I_2 = h$, we obtain, in terms of the invariants, the relation $3a_1 - a_2 = h$, which implies $a_2 = 3a_1 - h$. Substituting this expression into \eqref{hred31} and carrying out the algebraic simplifications, we arrive at
$$ {\mathscr H} =\frac{1}{2}\left[ (A+6B+9C)a_1^2-(2B+6C)ha_1\right]+\delta a_3+h+Ch^2.$$

Introducing the parameters $D = \tfrac{1}{2}(A+6B+9C)$ and $M = B+3C$, and dropping constant terms, the reduced averaged Hamiltonian takes the form
$$
\bar{R}=Da_1^2-Mha_1+\delta a_3.
$$
In KAM theory, the parameter $D$ is referred to as the \emph{twist coefficient}; see \cite{MeyerOffin}.

As in the case $k=2$, the system depends on the two masses $(m_1,m_2)$. 
The loci of parameter values for which the fundamental frequencies of the quadratic Hamiltonian part are rationally related form curves in the mass plane. 
Along each $3$:$-1$ resonance curve, the constant $\delta$ must be evaluated for each parameter pair.
For all pairs considered, we find $\delta \neq 0$; by continuity, this property extends along the entire curve.
Consequently, singular reduction theory is applicable for all admissible mass values.
It is therefore sufficient to select a representative point on the $3\!:\!-1$ resonance curve and carry out the analysis at that point. 
The resulting reduced dynamics is topologically equivalent along the curve, since the same geometric structure appears throughout.  

The reduced averaged system is defined on the orbit space $\mathbb{O}$, given by
\begin{equation}\label{orbspace31}
a_3^2 + a_4^2 = a_1(3a_1-h)^3, 
\qquad a_1 \geq 0, 
\qquad a_1 \geq \tfrac{1}{3}h.
\end{equation}
As before, the geometry of the problem is captured by the intersection of the orbit space surface (for different values of $h$) with the level surfaces of the reduced averaged Hamiltonian, 
$\bar{R}=\bar{h}$, for all $\bar{h} \in \mathbb{R}^3$.  

In the neighborhood of the origin in $(a_1,a_3,a_4)$ coordinates, the orbit space $\mathbb{O}$ exhibits distinct local structures: 
it is differentiable for $h<0$, a rotated parabola for $h=0$, and a rotated cusp for $h>0$. 
The surface $\bar{R}=\bar{h}=\mathrm{const}$ is simply a parabola translated along the $a_4$ direction.  

We now compute the vector field associated with $\bar{R}$ using the Poisson brackets in Table~\ref{tabrs1}:
\begin{align*}
\dot{a}_1 & =\{a_1,\bar{R}\}=\delta \{a_1,a_3\}=-2\delta a_4,\\
\dot{a}_3 & = \{a_3,\bar{R}\}=D\{a_3,a^2_1\}-Mh\{a_3,a_1\}\\
                & =2Da_1\{a_3,a_1\}-2Mha_4=4Da_1a_4-2Mha_4,\\
\dot{a}_4 & = \{a_4,\bar{R}\}=D\{a_4,a_1^2\}-Mh\{a_4,a_1\}+\delta \{a_4,a_3\}\\
               & =-4Da_1a_3+M ha_3-\delta a_2^2(9a_1+a_2).
\end{align*}

In the last expression, the dependence on $a_2$ can be removed, by using  $3a_1-a_2=h$. 
The resulting system reads
\begin{equation} \label{ham_sys4}
	\begin{split}
		\dot{a}_1& =-2\delta a_4,  \\
		\dot{a}_3&=2a_4(2Da_1-M h), \\
		\dot{a}_4&=-2a_3(2Da_1-Mh)-\delta (3a_1-h)^2(12a_1-h).
	\end{split}
\end{equation}
This system determines the precise flow on the orbit space. To locate equilibrium points, note that if $\delta \neq 0$, then necessarily $a_4=0$ so that the first two equations vanish, leaving the orbit space relation in $(a_1,a_3)$:
$$
a_3^2=a_1(3a_1-h)^3.
$$
From the third equation of \eqref{ham_sys4}, equilibria must satisfy
$$
2a_3(2Da_1-M h)+\delta (3a_1-h)^2(12a_1-h)=0.
$$
Solving for $a_3$ and substituting into the orbital constraint yields a polynomial condition in $a_1$:
$$\delta^2(3a_1-h)^4(12a_1-h)^2=4a_1(3a_1-h)^3(2Da_1-Mh)^2.$$
Simplification gives the cubic polynomial
\begin{equation}\label{eq_cubb}
	16(D^2-27\delta^2)a_1^3+8h(27\delta^2-2DM)a_1^2+h^2(4M^2-27\delta^2)a_1+\delta^2h^3=0,\end{equation}
whose coefficients depend on the parameters $D$, $\delta$, $M$, and $h$.  
Introducing scaled parameters $\alpha = D/\delta$ and $\beta = M/\delta$, and dividing by $\delta^2$, this reduces to
$$ 16(\alpha^2-27)a_1^3+8h(27-2\alpha\beta)a_1^2+h^2(4\beta^2-27)a_1+h^3=0.$$

We seek roots of the polynomial under the constraint $a_1 \geq \min \{0, h/3 \}$.  
The number of admissible roots changes when the polynomial and its derivative with respect to $a_1$ share a common factor, a condition that can be detected by computing their resultant.  
Applying the resultant method in $a_1$, we obtain the bifurcation condition in terms of $h$ and the system parameters:
$$
1024(\alpha^2-27)(\alpha-6\beta)(729+108\alpha^2+648\beta^2-48\beta^4+8\alpha\beta(-81+4\beta^2))h^6=0,$$
which characterizes the common zeros of the cubic polynomial~\eqref{eq_cubb} and its derivative.

\begin{figure}[H]\centering
	\subfigure[]{\includegraphics[width=0.4\textwidth]{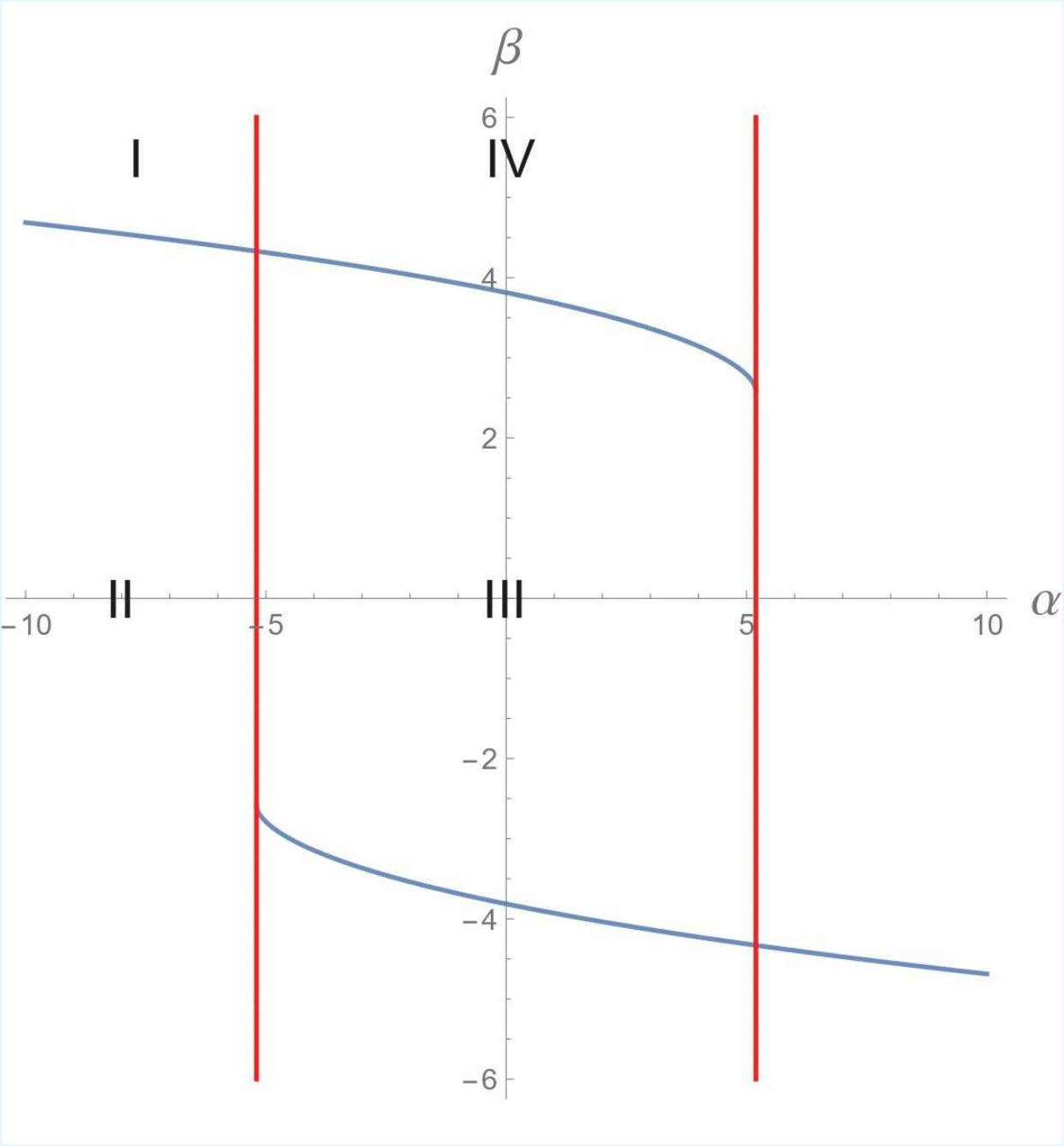}}\qquad 
	\subfigure[]{\includegraphics[width=0.4\textwidth]{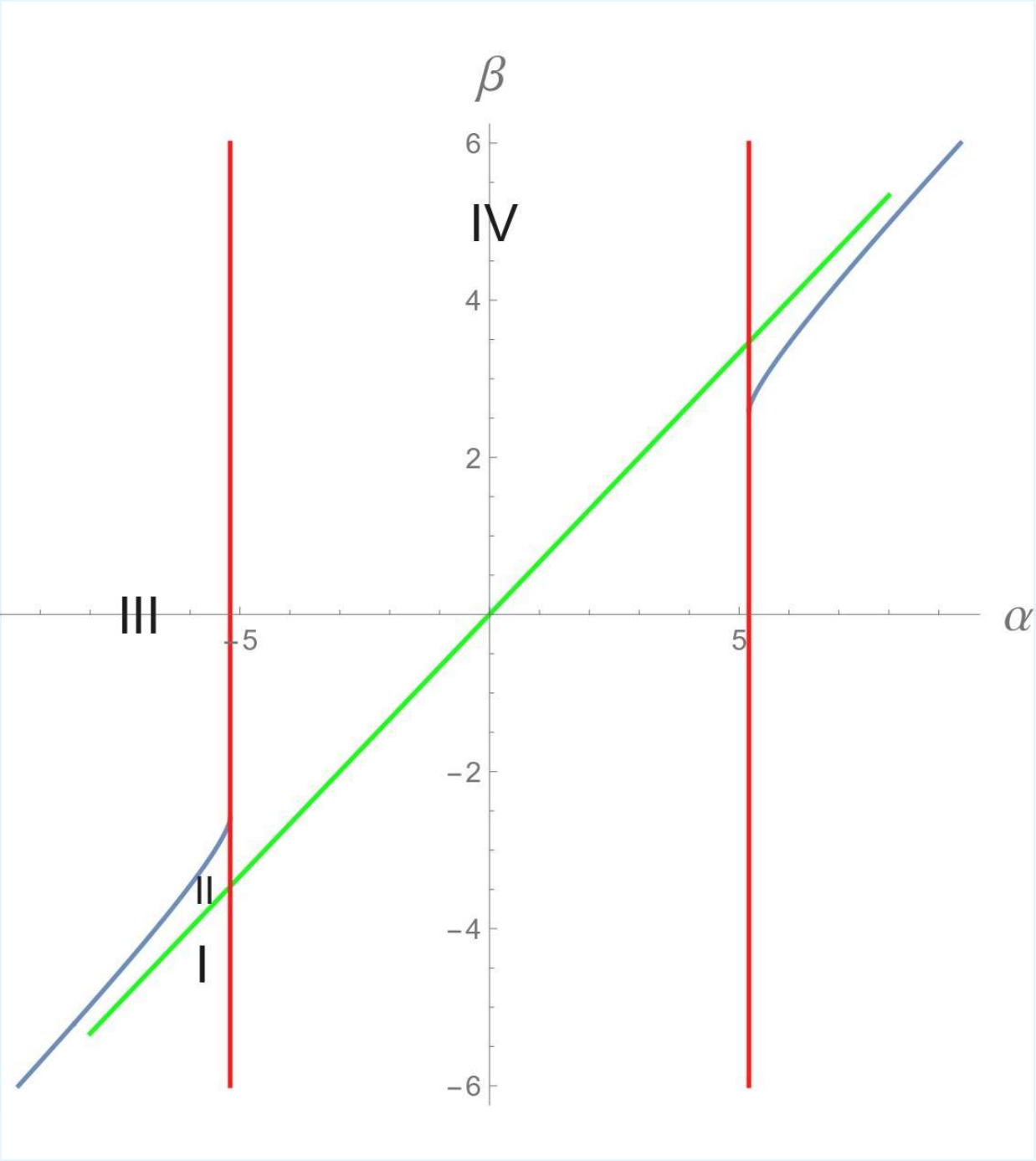}}	
	\caption{Bifurcation diagrams for: (a) $h<0$,  (b) $ h>0 $.} \label{rs3b}
\end{figure}

The number of critical points of system~\eqref{ham_sys4} in the orbit space changes whenever a root of the cubic polynomial~\eqref{eq_cubb} 
coincides with the maximum at $a_1 = h/3$ for $h \geq 0$, that is, when
$$\frac{4}{27}(2\alpha-3\beta)^2 h^3 = 0.$$
Accordingly, we distinguish three regimes, $h<0$, $h=0$, and $h>0$, and analyze the corresponding equilibria and their bifurcations as the parameters $\alpha$ and $\beta$ vary. The simplest case arises when $h=0$, for which equation~\eqref{eq_cubb} reduces to
$$16a_1^3(\alpha^2-27) = 0,$$
yielding either the trivial root $a_1 = 0$ or a bifurcation at $\alpha^2 = 27$.
For $h \neq 0$, two qualitatively different bifurcation diagrams emerge, illustrated in Figures~\ref{rs3b}(a) and~\ref{rs3b}(b) for $h < 0$ and $h > 0$, respectively. The corresponding bifurcation curves in the parameter plane $(\alpha,\beta)$ are given by
\begin{itemize}
\item[] $\Lambda_1$: $\alpha^2 - 27 = 0$ \; (red lines)
\item[] $\Lambda_2$:  $729 + 108\alpha^2 + 648\beta^2 - 48\beta^4 + 8\alpha\beta(-81 + 4\beta^2) = 0$
 \; (blue curves)
 \item[]  $\Lambda_3$: $2\alpha - 3\beta = 0$ \; (green line).
 \end{itemize}
These bifurcation diagrams are symmetric with respect to the origin in parameter space. The blue curves $\Lambda_2$ correspond to saddle-center bifurcations of critical points, leading to the emergence of periodic orbits. Along the red lines $\Lambda_1$, the cubic term in the governing polynomial vanishes, leaving only two real roots. The green line $\Lambda_3$ marks a \emph{peak bifurcation}, where a root coincides with the local maximum of the potential.

\begin{figure}[H]\centering
	\subfigure[]{\includegraphics[width=0.22\textwidth]{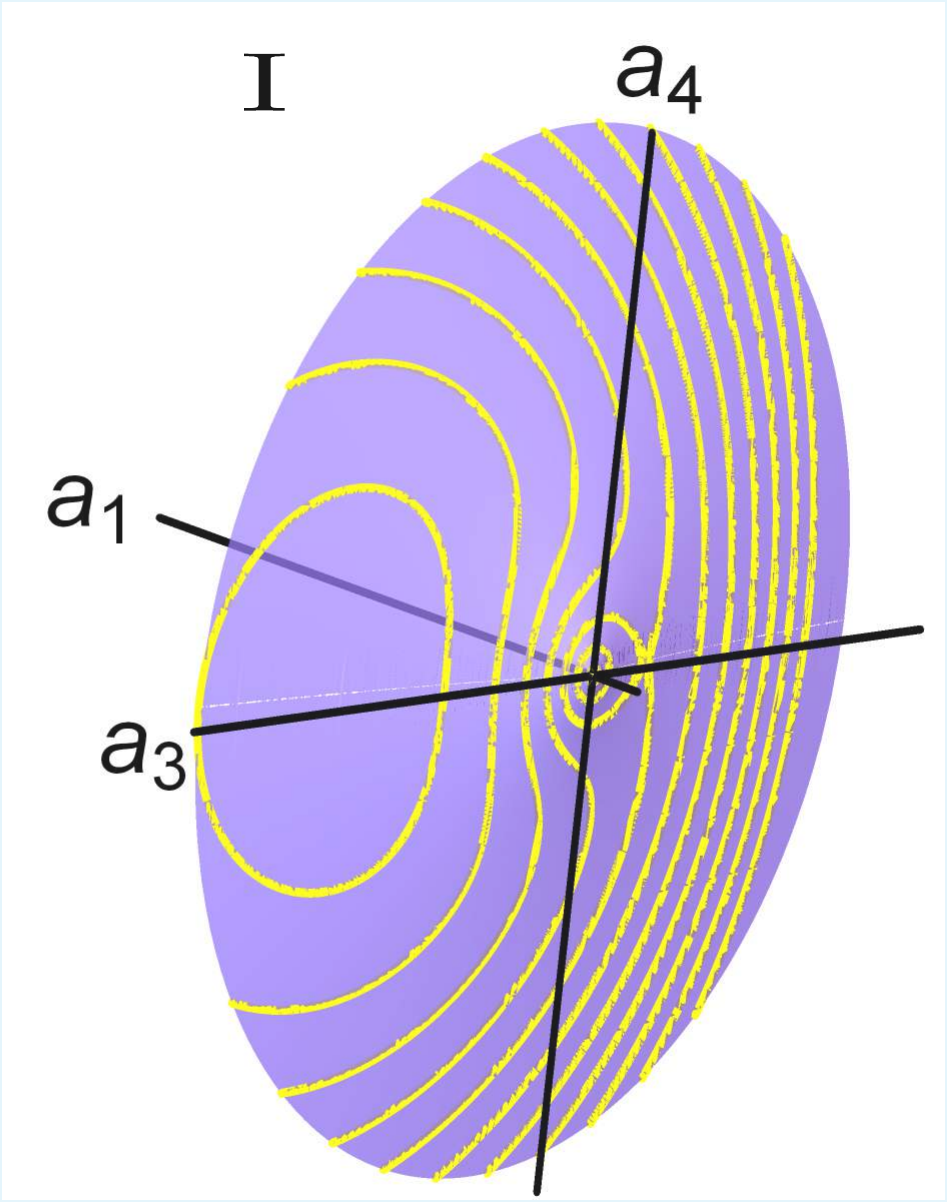}} \hspace{0.6cm}
	\subfigure[]{\includegraphics[width=0.22\textwidth]{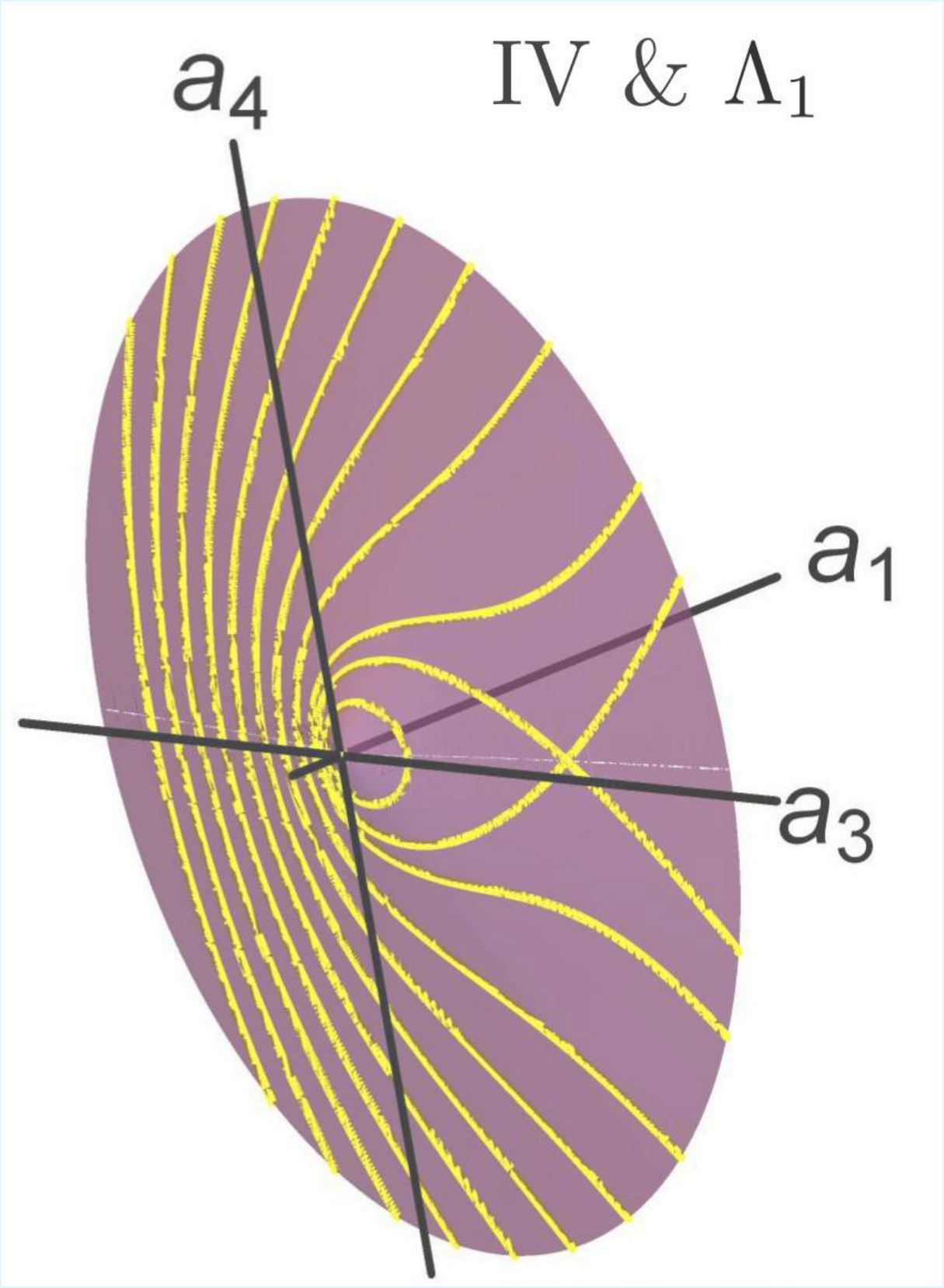}} \hspace{0.6cm}
	\subfigure[]{\includegraphics[width=0.22\textwidth]{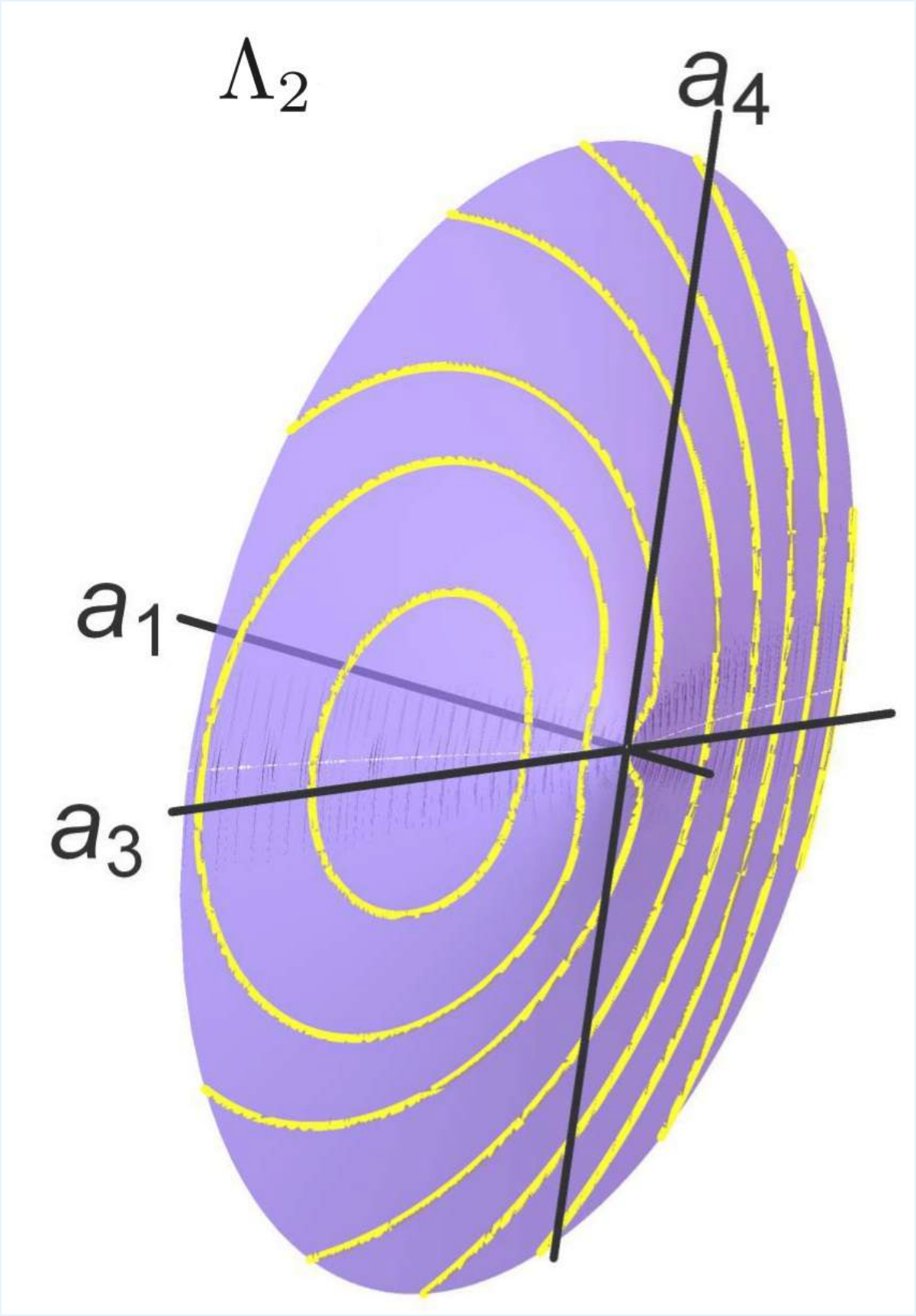}} \hspace{0.6cm}
	\subfigure[]{\includegraphics[width=0.22\textwidth]{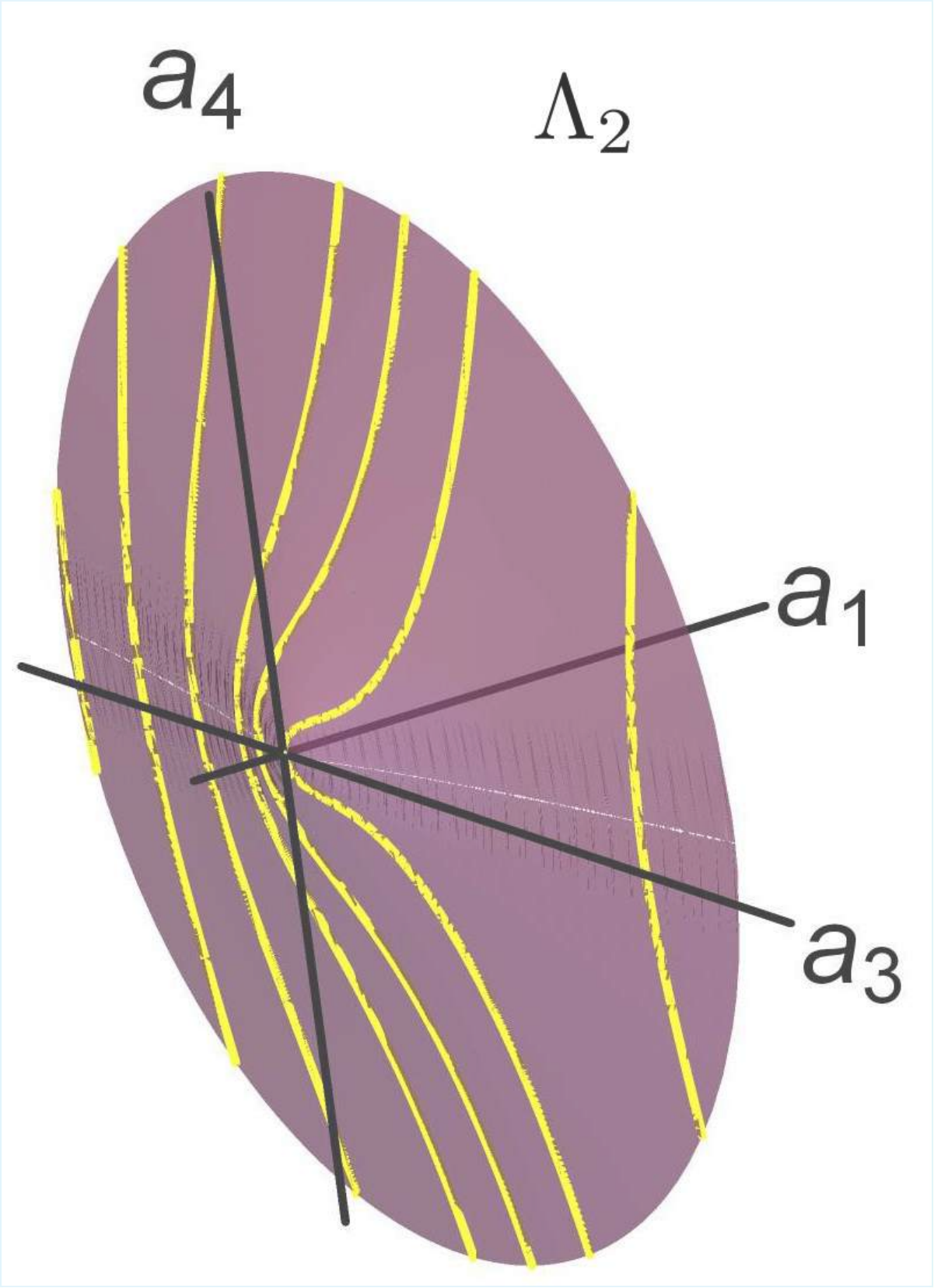}} \hspace{0.6cm}
	\subfigure[]{\includegraphics[width=0.22\textwidth]{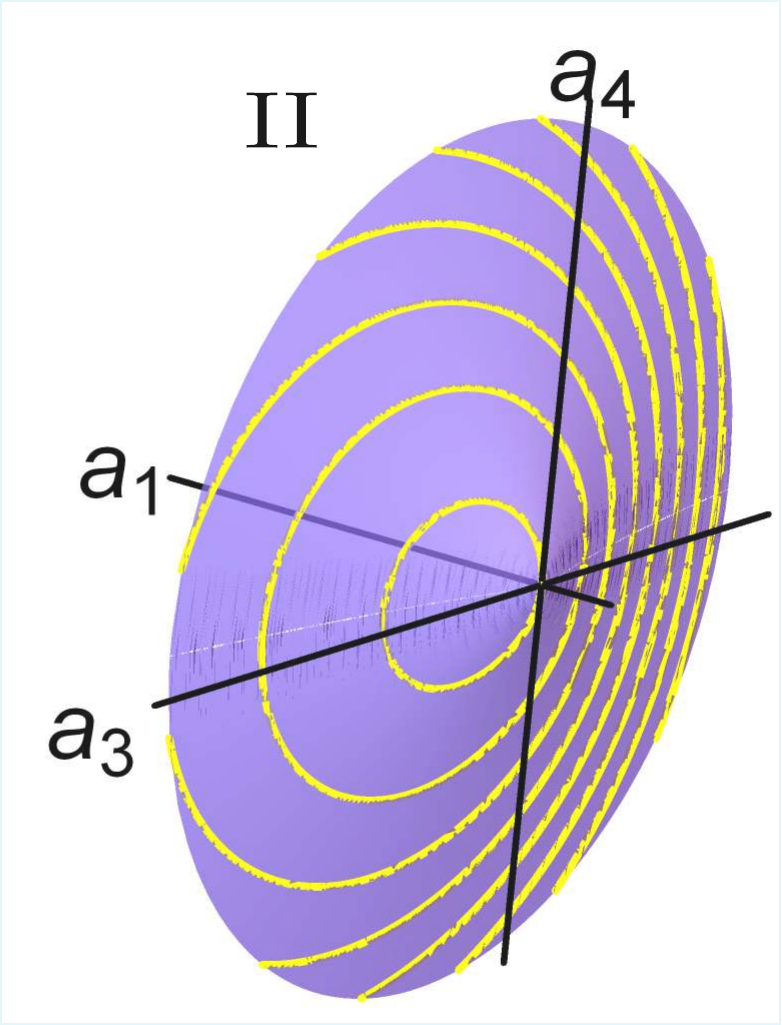}} \hspace{0.6cm}
	\subfigure[]{\includegraphics[width=0.19\textwidth]{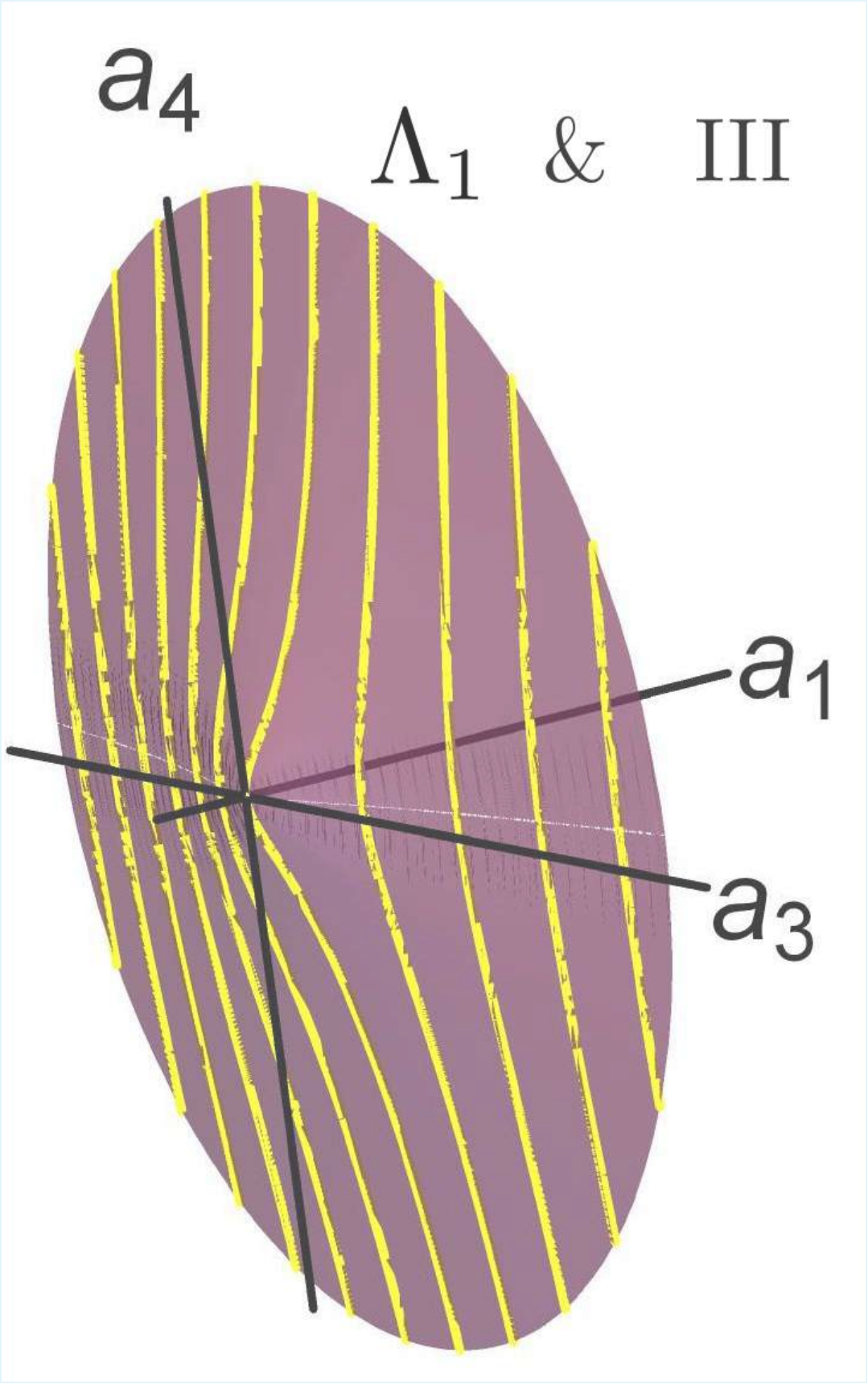}}
	\caption{Flow in the orbit space for the $3$:$-1$  resonance when  $h<0$.} \label{rs3d}
\end{figure}

In what follows, we describe the orbital flow for each of the three cases. We begin with the case $h < 0$, which is illustrated in Figure~\ref{rs3d}. In this scenario, the orbit space $\mathbb{O}$ is smooth, and the equilibria correspond to periodic orbits with period $2\pi$. We now turn to the analysis of the equilibria in each region of the bifurcation plane $(\alpha,\beta)$:
In region~I, equation~\eqref{eq_cubb} possesses three distinct positive roots, which correspond in the orbit space to two centers (stable periodic orbits) and one saddle (an unstable periodic orbit). According to KAM theory, the invariant curves surrounding the centers give rise to two-dimensional KAM tori.  
Along the blue curves ($\Lambda_2$), one of the centers collides with the saddle, producing a cusp-type degenerate equilibrium associated with a degenerate periodic orbit, while the other center remains enclosed by KAM tori. Crossing into region~II, this degenerate point disappears and only the surviving center persists.  
In region~IV, which also includes the red bifurcation lines ($\Lambda_1$), the orbit space contains precisely one center and one saddle. Once again, on the blue curves ($\Lambda_2$), these two equilibria coalesce into a cusp-type degenerate point.  
Finally, in region~III the situation is simpler: no equilibria are present.  

\begin{figure}[H]\centering
	\subfigure[]{\includegraphics[width=0.22\textwidth]{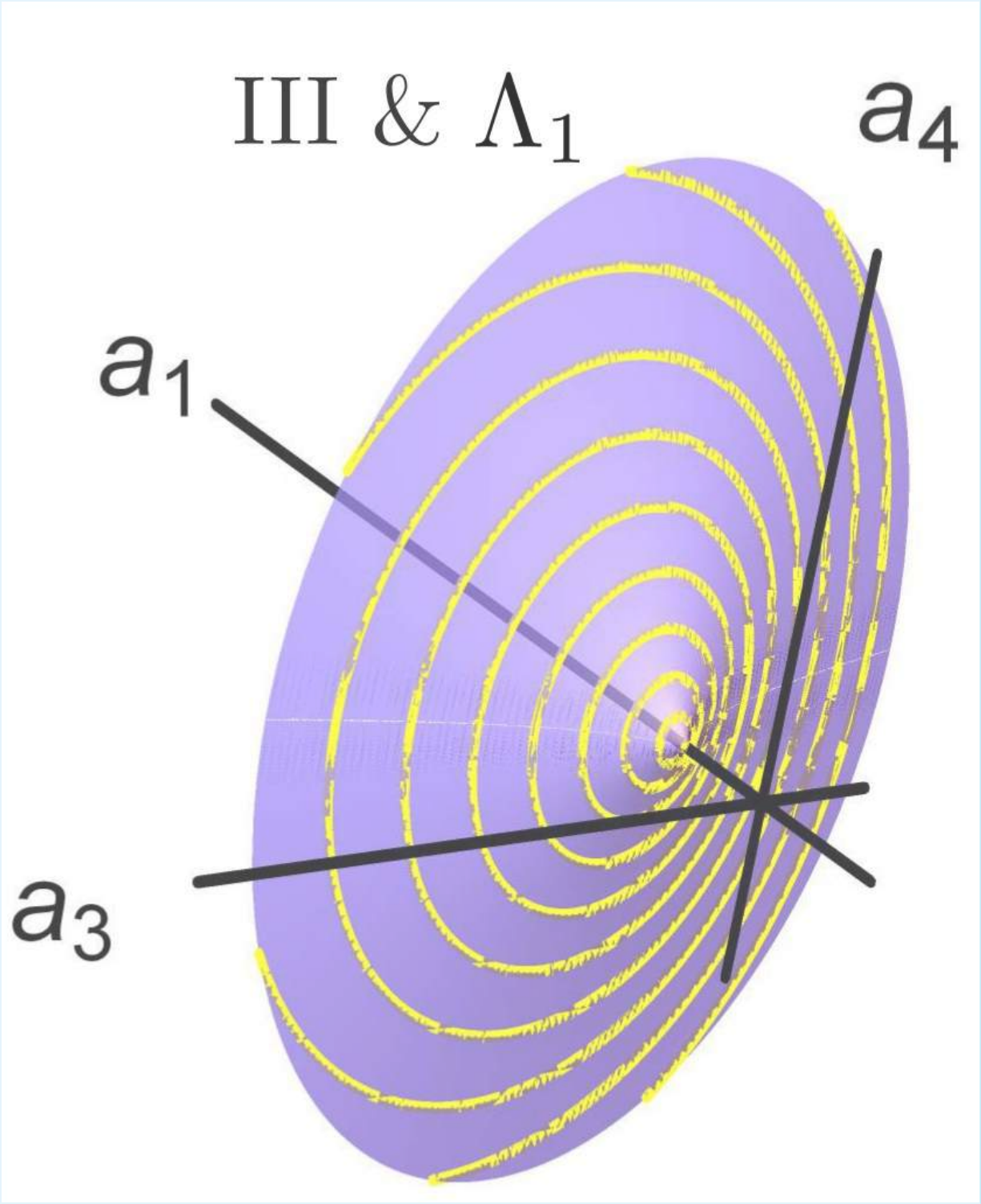}} \hspace{0.6cm} 
	\subfigure[]{\includegraphics[width=0.22\textwidth]{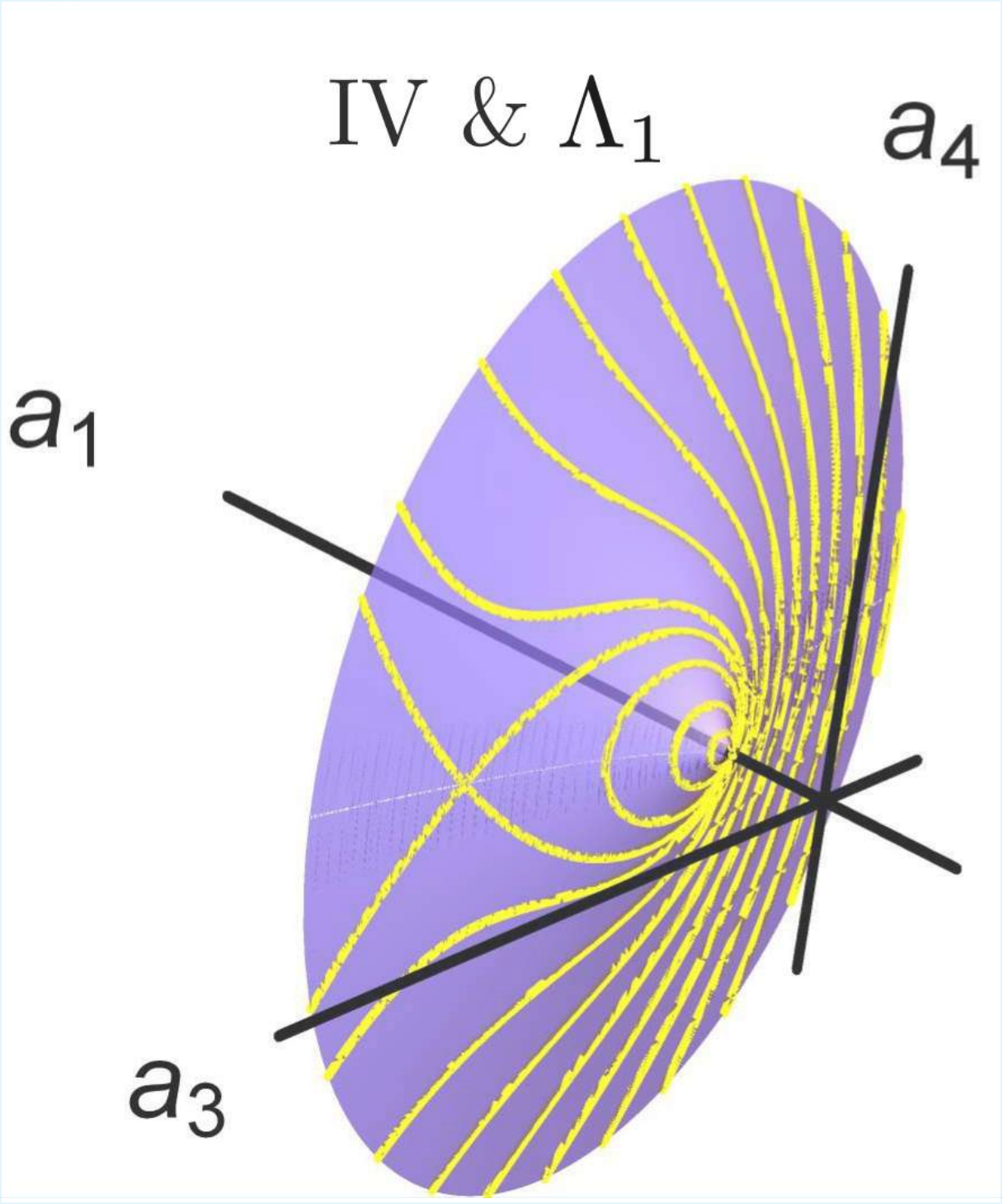}}	 \hspace{0.6cm} 
	\subfigure[]{\includegraphics[width=0.22\textwidth]{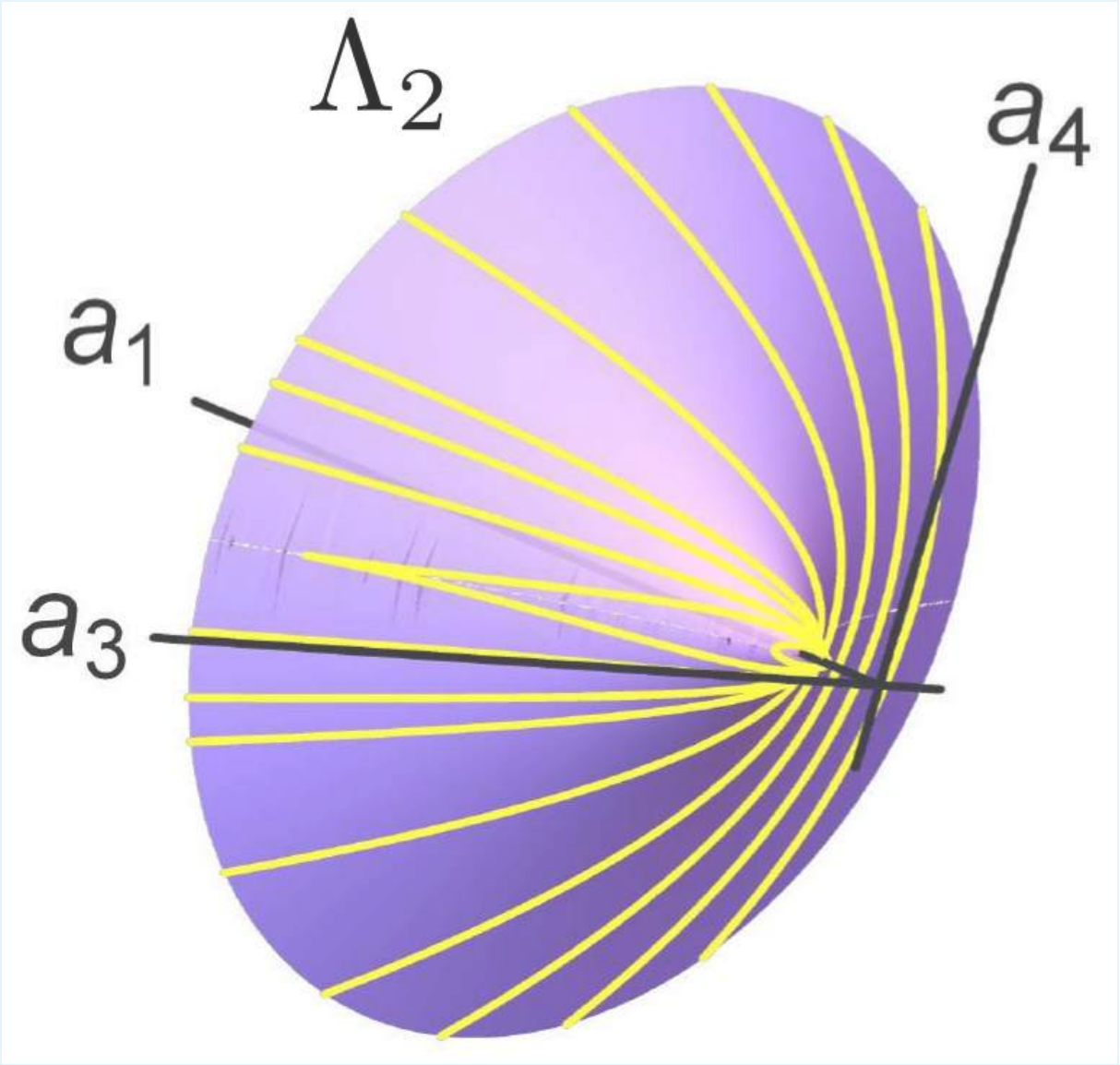}} \hspace{0.7cm}
	\subfigure[]{\includegraphics[width=0.22\textwidth]{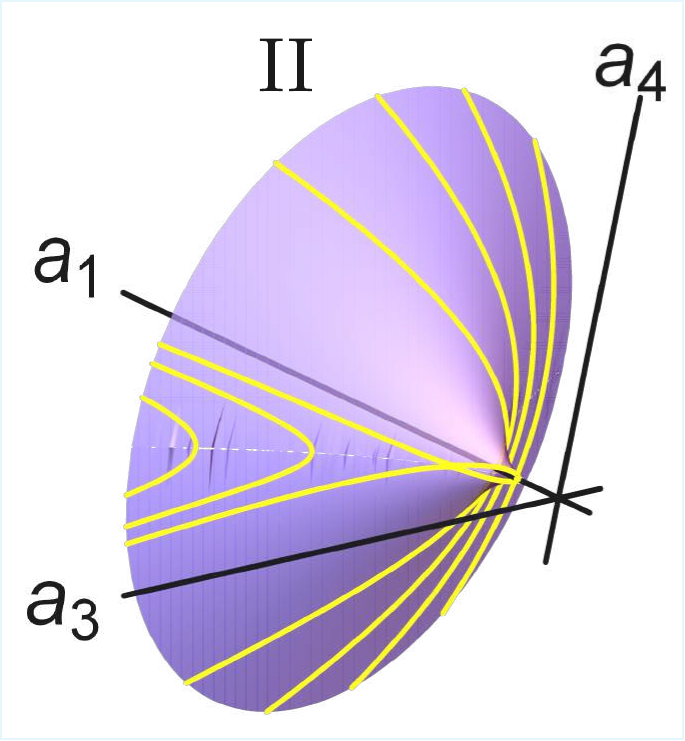}} \hspace{0.6cm}
	\subfigure[]{\includegraphics[width=0.22\textwidth]{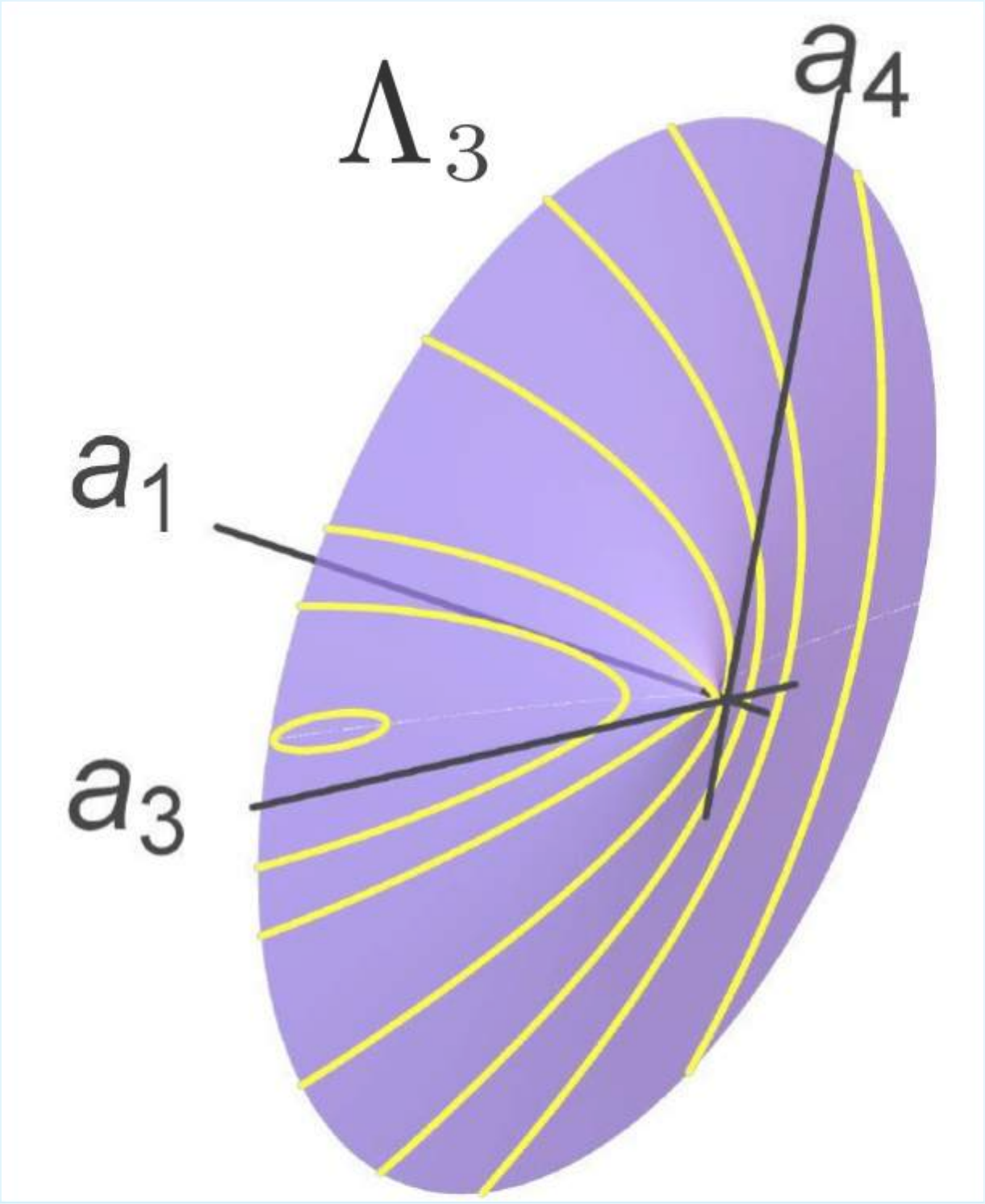}} \hspace{0.6cm} 
	\subfigure[]{\includegraphics[width=0.22\textwidth]{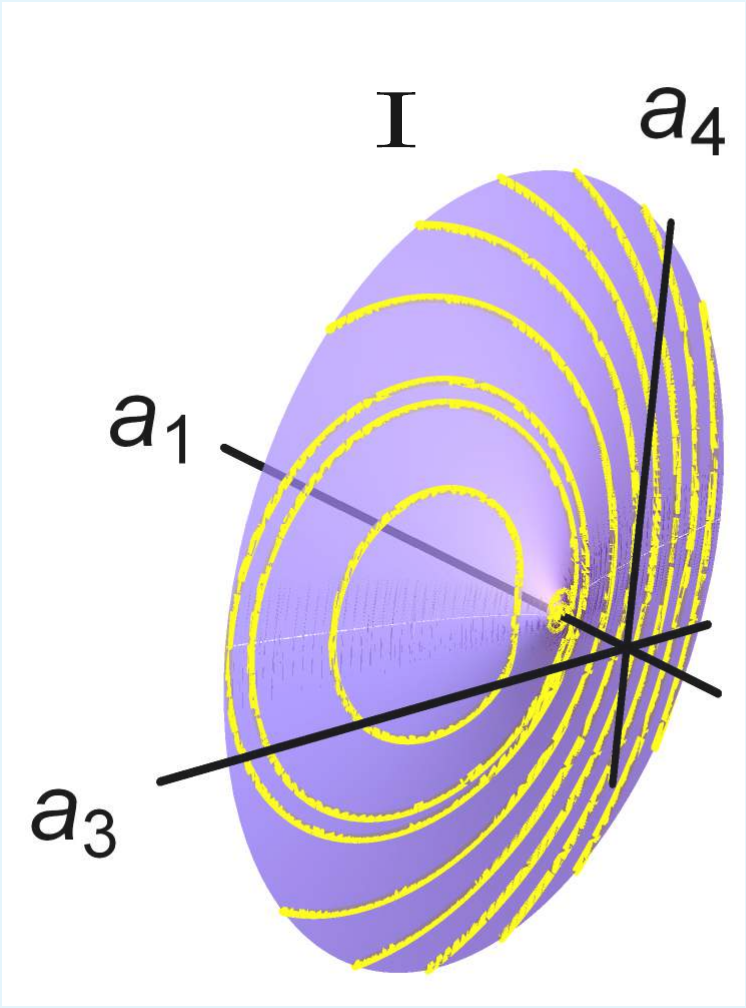}}
	\caption{Flow in the orbit space for the $3$:$-1$  resonance when $h>0$.} \label{rs3c}
\end{figure}

For $h>0$, as illustrated in Figure~\ref{rs3c}, the point $a_1 = h/3$ always defines an equilibrium.  
The key question is whether the cubic admits additional roots greater than $h/3$.  
In Region~I, two such roots are present, leading to a total of three equilibria: the peak, which acts as a center; a saddle situated nearby; and another center located at a greater distance.
Along the green curve, the saddle merges with the peak, so that the configuration of Region~I collapses into that of Region~II.  
Moving instead to the blue curve, the saddle coalesces with the outer center, and both equilibria vanish, leaving in Region~III only the peak as a surviving center, up to the red boundary.  Beyond this red line, in Region~IV, a new saddle appears while the peak retains its role as a center.

\section{Conclusions}\label{sec:conclusions}
This work establishes, for the first time, a complete nonlinear stability classification of
the elliptic equilibria $L_3$, $L_5$, and $L_6$ across the full admissible mass plane of the
unequal-mass equilateral restricted four-body problem. Three results summarize the main
contributions. First, nonresonant stability is confirmed throughout the admissible mass plane
via the quartic Arnold determinant $D_4$, extended to the sixth-order determinant $D_6$ where
$D_4$ degenerates. Second, both the $2$:$-1$ and $3$:$-1$ resonances are shown to be generically
unstable, with the single exception of a stable segment of $L_5$ under the $3$:$-1$ resonance;
this is the first explicit treatment of the $3:-1$ resonance in the unequal-mass problem.
Third, singular reduction of the resonant normal forms shows that the associated orbit spaces
share a universal geometric structure with the classical restricted three-body problem, with
the mass parameters entering only through the bifurcation thresholds $\Lambda_1$, $\Lambda_2$,
and $\Lambda_3$ of the reduced Hamiltonians; this confirms that the local instability detected
in Section~\ref{sec_nonstab} corresponds to genuine global bifurcation phenomena, governed by the appearance,
coalescence, and disappearance of centers and saddles associated with families of periodic
orbits, rather than to an artifact of the local normal-form analysis. As noted in
Section~\ref{sec0}, the present study is restricted to the planar problem and to the two lowest-order
resonances; extending this framework to the spatial ERFBP and to higher-order resonances,
where richer resonance combinations and additional degrees of freedom arise, is left for
future work.

\section*{Acknowledgements}  M. Alvarez-Ramírez gratefully acknowledges the support of Proyectos CBI UAM-I 2026.
\appendix
\section{Stability theorems}\label{ap1}
The versions used here are adapted from \cite{MeyerOffin}.

\begin{theorem}[Arnold’s stability theorem] \label{teo_arnold} Consider a Hamiltonian $H$ with two degrees of freedom, expressed  in  action-angle variables $(I_1,I_2,\phi_1,\phi_1)$ as 
	\begin{equation}\label{H29}
		H=H_2+H_4 + \cdots + H_{2n}+H^*,
	\end{equation}
	where 
	\begin{enumerate}
		\item[{1.}] $H$ is real analytic in a neighborhood of the origin in $\mathbb{R}^4$.
		\item[{2.}]  $H_{2k}$,  $k= 1,\dots, n$  are homogeneous polynomial of degree $k$ in $I_1$, $I_2$, with real 
		coefficients. In particular
		\begin{enumerate}
			\item[{}] $H_2 = \omega_1 I_1-\omega_2 I_2$ with $\omega_1$, $\omega_2$ nonzero constant,
			\item[{}] $H_4=\frac{1}{2} (A I_1^2+2BI_1I_2 + CI_2^2)$, with $A,B,C$ constant.
		\end{enumerate}
		\item[{3.}]  $H^*$ has a series expansion that starts with terms at least of degree $2n+1$.
	\end{enumerate}
	Under these assumptions,  the origin  is  stable  for the system whose Hamiltonian is \eqref{H29}, provided that for some $k$, $1\leq k\leq n$, $D_{2k}= H_{2k}(\omega_2,\omega_1)\neq 0$ or, equivalently, provided that $H_2$ 
	does not divide $H_{2k}$. In particular, the equilibrium is stable if
	\begin{equation}\label{arnold_cond}
		D_4= \frac{1}{2} \{A \omega_2^2 + 2 B\omega_1\omega_2 + C\omega_1^2\}\neq 0.
	\end{equation}
	Moreover, arbitrarily close to the origin in $\mathbb{R}^4$, there are invariant tori and the flow on these invariant tori is the linear flow with irrational slope.
\end{theorem} 
However, there are  parameter values at which resonances occur, so that the Theorem \ref{teo_arnold} does not apply. 

Consider the case where the linear system is in $1$:$2$ resonance, that is, when the
linearized system has eigenvalues  $\pm i \omega_1$ and $\pm i \omega_2$ with $\omega_1 = 2\omega_2$. Let $\omega = \omega_2$. 
Assume the system has been normalized up to third-degree terms, meaning the Hamiltonian takes the following form:
\begin{equation}\label{ap_res_21}
H  = 2 \omega I_1 -\omega I_2 + \delta I_1^{1/2}I_2\cos\psi + H^\dagger,
 \end{equation}
where $ \psi=\phi_1+2\phi_2 $  and $H^\dagger (I_1,I_2,\phi_1,\phi_2)= \mathcal{O} ((I_1+I_2)^2)$.

\begin{theorem}[Alfriend-Markeev-Meyer]\label{teo_AMM}
If in the presence of $ 2 $:$-1$ resonance, the Hamiltonian system is in the normal form \eqref{ap_res_21}
with $ \delta\neq 0$ then the equilibrium is unstable. In fact, there is a neighborhood $O$ of the equilibrium such that any solution starting in $O$ and not on the Lyapunov center leaves $O$ in either positive or negative time. In particular, the small periodic solutions of the short period family are unstable.
\end{theorem}

Now, consider the system in the case when the linear system is in $1$:$3$ resonance, that is, when $\omega_1 = 3\omega_2$. 
Let $\omega = \omega_2$. Assuming that the system has been normalized up to terms of degree four, the Hamiltonian takes the following form:
\begin{equation}\label{ap_res_31}
H =  3\omega I_1- \omega I_2  +\delta I_1^{1/2} I_2^{3/2}\cos\psi +
 \frac{1}{2} \left\{  A I_1^2 + 2 B I_1I_2 + C I_2^2\right\} + H^\dagger, 
\end{equation}
where $\psi = \phi_1+3\phi_2$, $H^\dagger = {\mathcal O} ((I_1+I_2)^{5/2}))$. Let $D= A+6B+9C$, 
which is related to the Arnold determinant by $D_4= \frac{1}{2}D\omega^2$.

\begin{theorem}[Alfriend-Markeev]\label{teo31}
If in the presence of $3$:$-1$  resonance, the Hamiltonian system is in the normal form \eqref{ap_res_31}
 and if $ 6\sqrt{3}|\delta|>|D| $ then the equilibrium is unstable, whereas, if $ 6\sqrt{3}|\delta|<|D| $ then the equilibrium is stable.
\end{theorem}

\bibliographystyle{abbrvnat}

\begin{thebibliography}{99}

\bibitem{AlvarezBarrabes2014}
M. Alvarez-Ram{\'i}rez and E. Barrab{\'e}s,
Transport orbits in an equilateral restricted four-body problem,
\textit{Celest. Mech. Dyn. Astron.} \textbf{120}(4) (2014) 389--405.
doi: 10.1007/s10569-014-9594-z.

\bibitem{alvastuchi}
M. Alvarez-Ram{\'i}rez, J. E. F. Skea, and T. J. Stuchi,
Nonlinear stability analysis in an equilateral restricted four-body problem,
\textit{Astrophys. Space Sci.} \textbf{358} (2015) 1--11.
doi: 10.1007/s10509-015-2333-4.

\bibitem{amiet2002commensurate}
J.-P. Amiet and S. Weigert,
Commensurate harmonic oscillators: Classical symmetries,
\textit{J. Math. Phys.} \textbf{43}(8) (2002) 4110--4126.
doi: 10.1063/1.1488672.

\bibitem{arens}
R. F. Arenstorf,
Central configurations of four bodies with one inferior mass,
\textit{Celest. Mech.} \textbf{28}(1--2) (1982) 9--15.
doi: 10.1007/BF01230655.

\bibitem{arms1991universal}
J. M. Arms, R. H. Cushman, and M. J. Gotay,
A universal reduction procedure for Hamiltonian group actions,
in: J. E. Marsden and T. S. Ratiu (Eds.), \textit{The Geometry of Hamiltonian Systems}, pp. 33--51, Springer, New York, 1991.
doi: 10.1007/978-1-4613-9725-0\_4.

\bibitem{arnold1961}
V. I. Arnol'd,
The stability of the equilibrium position of a Hamiltonian system of ordinary differential equations in the general elliptic case,
\textit{Dokl. Akad. Nauk SSSR} \textbf{137}(2) (1961) 255--257.

\bibitem{Arnold2006}
V. I. Arnol'd, V. V. Kozlov, and A. I. Neishtadt,
Mathematical Aspects of Classical and Celestial Mechanics,
3rd ed., Encyclopaedia Math. Sci., Vol. 3, Springer, Berlin, 2006.

\bibitem{Balta1}
A. N. Baltagiannis and K. E. Papadakis,
Equilibrium points and their stability in the restricted four-body problem,
\textit{Int. J. Bifurcation Chaos} \textbf{21}(8) (2011) 2179--2193.
doi: 10.1142/S0218127411029707.

\bibitem{Bardin2018}
B. S. Bardin and P. A. Esipov,
Investigation of Lyapunov stability of a central configuration in the restricted four-body problem,
in: \textit{AIP Conf. Proc.}, Vol. 1959(1), pp. 040004, 2018.
doi: 10.1063/1.5034607.

\bibitem{BardinVolkov2024}
B. S. Bardin and E. V. Volkov,
The Lyapunov stability of central configurations of the planar circular restricted four-body problem,
\textit{Cosmic Res.} \textbf{62} (2024) 388--400.
doi: 10.1134/S0010952524600677.

\bibitem{BudProk2011}
D. A. Budzko and A. N. Prokopenya,
On the stability of equilibrium positions in the circular restricted four-body problem,
in: V. Gerdt, E. W. Mayr, and E. V. Vorozhtsov (Eds.), \textit{Computer Algebra in Scientific Computing}, pp. 88--100, Springer, Berlin, 2011.
doi: 10.1007/978-3-642-23568-9\_8.

\bibitem{cabral}
H. Cabral and K. Meyer,
Stability of equilibria and fixed points of conservative systems,
\textit{Nonlinearity} \textbf{12}(5) (1999) 1351--1362.
doi: 10.1088/0951-7715/12/5/309.

\bibitem{cushman1997global}
R. H. Cushman and L. M. Bates,
Global Aspects of Classical Integrable Systems,
Birkh{\"a}user, Basel, 1997.
doi: 10.1007/978-3-0348-0918-4.

\bibitem{Dirichlet}
G. L. Dirichlet,
{\"U}ber die Stabilit{\"a}t des Gleichgewichts,
\textit{J. Reine Angew. Math.} \textbf{32} (1846) 85--88.

\bibitem{FerrerPalacian2026}
S. Ferrer-Bened{\'i} and J. F. Palaci{\'a}n,
Singular Reduction of Resonant Hamiltonians with Two or More Degrees of Freedom,
\textit{Qual. Theory Dyn. Syst.} \textbf{25} (2026) 136.
doi: 10.1007/s12346-026-01560-7.

\bibitem{Gasch}
M. Gascheau,
Examen d'une classe d'{\'e}quations diff{\'e}rentielles et application {\`a} un cas particulier du probl{\`e}me des trois corps,
\textit{C. R. Acad. Sci. Paris} \textbf{16}(7) (1843) 393--394.

\bibitem{KepleyMireles2019}
S. Kepley and J. D. Mireles James,
Chaotic motions in the restricted four body problem via Devaney's saddle-focus homoclinic tangle theorem,
\textit{J. Differ. Equ.} \textbf{266}(4) (2019) 1709--1755.
doi: 10.1016/j.jde.2018.08.007.

\bibitem{leandro}
E. S. G. Leandro,
On the central configurations of the planar restricted four-body problem,
\textit{J. Differ. Equ.} \textbf{226} (2006) 323--351.
doi: 10.1016/j.jde.2005.10.015.

\bibitem{Markeev}
A. P. Markeev,
Tochki libratsii v nebesnoi mekhanike i kosmodinamike (Libration Points in Celestial Mechanics and Space Dynamics),
Nauka, Moscow, 1978, In Russian.

\bibitem{MeyerOffin}
K. R. Meyer and D. C. Offin,
Introduction to Hamiltonian Dynamical Systems and the $N$-Body Problem,
Appl. Math. Sci., Vol. 90, Springer, Cham, 2017.

\bibitem{palacian2018}
K. R. Meyer, J. Palaci{\'a}n, and P. Yanguas,
Singular reduction of resonant Hamiltonians,
\textit{Nonlinearity} \textbf{31}(6) (2018) 2854--2894.
doi: 10.1088/1361-6544/aab591.

\bibitem{Moulton1900}
F. R. Moulton,
On a class of particular solutions of the problem of four  bodies,
\textit{Trans. Amer. Math. Soc.} \textbf{1}(1) (1900) 17--29.
doi: 10.2307/1986399.

\bibitem{olver}
P. J. Olver,
Classical Invariant Theory,
London Math. Soc. Student Texts, Vol. 44, Cambridge Univ. Press, Cambridge, 1999.
doi: 10.1017/CBO9780511623660.

\bibitem{palacian2015}
J. F. Palaci{\'a}n, F. Sayas, and P. Yanguas,
Flow reconstruction and invariant tori in the spatial three-body problem,
\textit{J. Differ. Equ.} \textbf{258}(6) (2015) 2114--2159.
doi: 10.1016/j.jde.2014.12.001.

\bibitem{pede}
P. Pedersen,
Librationspunkte im restringierten Vierk{\"o}rperproblem,
\textit{Publ. Mindre Medd. K{\o}benhavns Obs.} \textbf{137} (1944) 1--80.

\bibitem{Routh}
E. J. Routh,
On Laplace's three particles, with a supplement on the stability of steady motion,
\textit{Proc. Lond. Math. Soc.} \textbf{1}(1) (1874) 86--97.
doi: 10.1112/plms/s1-6.1.86.

\bibitem{simo}
C. Sim{\'o},
Relative equilibrium solutions in the four body problem,
\textit{Celest. Mech.} \textbf{18}(2) (1978) 165--184.
doi: 10.1007/BF01228714.

\bibitem{sjamaar1991stratified}
R. Sjamaar and E. Lerman,
Stratified symplectic spaces and reduction,
\textit{Ann. Math.} \textbf{134}(2) (1991) 375--422.
doi: 10.2307/2944350.

\bibitem{zepAl}
J. A. Zepeda Ram{\'i}rez, M. Alvarez-Ram{\'i}rez, and A. Garc{\'i}a,
Nonlinear stability of equilibrium points in the planar equilateral restricted mass-unequal four-body problem,
\textit{Internat. J. Bifur. Chaos Appl. Sci. Engrg.} \textbf{31}(11) (2021) 2130031.
doi: 10.1142/S0218127421300317.

\bibitem{zepeAl2}
J. A. Zepeda Ram{\'i}rez, M. Alvarez-Ram{\'i}rez, and A. Garc{\'i}a,
A note on the nonlinear stability of equilibrium points in the planar equilateral restricted mass-unequal four-body problem,
\textit{Internat. J. Bifur. Chaos Appl. Sci. Engrg.} \textbf{32}(2) (2022) 2250029.
doi: 10.1142/S0218127422500298.

\bibitem{zepeAl3}
J. A. Zepeda Ram{\'i}rez and M. Alvarez-Ram{\'i}rez,
Equilibrium points and their linear stability in the planar equilateral restricted four-body problem: A review and new results,
\textit{Astrophys. Space Sci.} \textbf{367}(8) (2022) 77.
doi: 10.1007/s10509-022-04108-8.

\bibitem{zotos}
E. E. Zotos,
Exploring the location and linear stability of the equilibrium points in the equilateral restricted four-body problem,
\textit{Internat. J. Bifur. Chaos Appl. Sci. Engrg.} \textbf{30}(10) (2020) 2050155.
doi: 10.1142/S0218127420501552.

\end{thebibliography}

\end{document}